\documentclass[a4paper,11pt]{article}

\usepackage{geometry}
\usepackage{amsfonts}
\usepackage{amsthm}
\usepackage{amssymb}
\usepackage{amsmath}
\usepackage{upref}
\usepackage{mathrsfs}
\usepackage{color}
\usepackage[citecolor=blue,colorlinks=true]{hyperref}
\usepackage{enumitem}
\usepackage{graphicx}
\usepackage{tikz}

\theoremstyle{plain}
\newtheorem{theorem}{Theorem}[section]
\newtheorem{proposition}[theorem]{Proposition}
\newtheorem{lemma}[theorem]{Lemma}
\newtheorem{corollary}[theorem]{Corollary}

\theoremstyle{definition}
\newtheorem{definition}{Definition}[section]

\theoremstyle{remark}
\newtheorem{remark}{Remark}[section]

\numberwithin{equation}{section}
\def\N{\mathbb{N}}
\def\Z{\mathbb{Z}}

\def\R{\mathbb{R}}

\def\ds{\displaystyle} 
\def\div{{\rm div}}
\def\ocirc#1{\ifmmode\setbox0=\hbox{$#1$}\dimen0=\ht0
    \advance\dimen0 by1pt\rlap{\hbox to\wd0{\hss\raise\dimen0
    \hbox{\hskip.2em$\scriptscriptstyle\circ$}\hss}}#1\else
    {\accent"17 #1}\fi} 

\def\eps{\varepsilon}
\def\<{\langle}
\def\>{\rangle}
\def\F{\mathcal{F}}
\def\GG{\mathbf{G}}
\def\P{\mathbb{P}}
\def\E{\mathbb{E}}
\def\T{\mathbb{T}}

\def\order{\mathrm{d}}

\begin{document}

\title{Convergence of the Finite Volume Method for scalar conservation laws with multiplicative noise: an approach by kinetic formulation}
\author{Sylvain Dotti\thanks{Aix-Marseille Universit\'e, CNRS, Centrale Marseille, I2M, UMR 7373, 13453 Marseille France }{ } and Julien Vovelle\thanks{Univ Lyon, Universit\'e Claude Bernard Lyon 1, CNRS UMR 5208, Institut Camille Jordan, 43 blvd. du 11 novembre 1918, F-69622 Villeurbanne cedex, France. Julien Vovelle was supported by ANR projects STOSYMAP and STAB.}}
\maketitle

\begin{abstract} We prove the convergence of the explicit-in-time Finite Volume method with monotone fluxes for the approximation of scalar first-order conservation laws with multiplicative, compactly supported noise. 
\end{abstract}

{\bf Keywords:}  Finite Volume Method, stochastic scalar conservation law, kinetic formulation

{\bf MSC Number:} 65M08 (35L60 35L65 35R60 60H15 65M12)

\tableofcontents

\section{Introduction}\label{secIntro}

\textbf{Stochastic first-order scalar conservation law.} Let $(\Omega,\F,\P,(\F_t),(\beta_k(t)))$ be a stochastic basis and let $T>0$. Consider the first-order scalar conservation law with stochastic forcing
\begin{equation}\label{stoSCL}
du(x,t)+\div(A(u(x,t)))dt=\Phi(x,u(x,t)) dW(t),\quad x\in\T^N, t\in(0,T).
\end{equation}
Equation~\eqref{stoSCL} is periodic in the space variable:  $x\in\T^N$ where $\T^N$ is the $N$-dimensional torus. The flux function $A$ in \eqref{stoSCL} is supposed to be of class $C^2$: $A\in C^2(\R;\R^N)$. We assume that $A$ and its derivatives have at most polynomial growth. Without loss of generality, we will assume also that $A(0)=0$. The right-hand side of \eqref{stoSCL} is a stochastic increment in infinite dimension. It is defined as follows (see \cite{DaPratoZabczyk92} for the general theory): $W$ is a cylindrical Wiener process, $W=\sum_{k\geq 1}\beta_k e_k$, where the coefficients $\beta_k$ are independent Brownian processes and $(e_k)_{k\geq 1}$ is a complete orthonormal system in a Hilbert space $H$. For each $x\in\T^N$, $u\in\R$, $\Phi(x,u)\in L_2(H,\R)$ is defined by $\Phi(x,u)e_k=g_k(x,u)$ where $g_k(\cdot,u)$ is a regular function on $\T^N$. Here, $L_2(H,K)$ denotes the set of Hilbert-Schmidt operator from the Hilbert space $H$ to an other Hilbert space $K$. Since $K=\R$ in our case, this set is isomorphic to $H$, thus we may also define
$$
\Phi(x,u)=\sum_{k\geq 1}g_k(x,u) e_k,
$$ 
the action of $\Phi(x,u)$ on $e\in H$ being given by $\<\Phi(x,u),e\>_H$. We assume $g_k\in C(\T^N\times\R)$, with the bounds 
\begin{align}
\GG^2(x,u)=\|\Phi(x,u)\|_H^2=\sum_{k\geq 1}|g_k(x,u)|^2\leq D_0(1+|u|^2),\label{D0}\\
\|\Phi(x,u)-\Phi(y,v)\|_H^2=\sum_{k\geq 1}|g_k(x,u)-g_k(y,v)|^2\leq D_1(|x-y|^2+|u-v|h(|u-v|)),\label{D1}
\end{align}
where $x,y\in\T^N$, $u,v\in\R$, and $h$ is a continuous non-decreasing function on $\R_+$ such that $h(0)=0$. We assume also $0\leq h(z)\leq 1$ for all $z\in\R_+$. \bigskip

\textit{Notation:} in what follows, we will use the convention of summation over repeated indices $k$. For example, we write $W=\beta_k e_k$.
\medskip

\textbf{Compactly supported multiplicative noise.} In this paper, we study the numerical approximation of \eqref{stoSCL}: our aim is to prove the convergence of the Finite Volume method with monotone fluxes, see Theorem~\ref{th:mainthm}. Our analysis will be restricted to the case of \textit{multiplicative noise with compact support}. Indeed, from Section~\ref{secFVscheme} to Section~\ref{sec:Convergence}, we will work under the following hypothesis: there exists $a,b\in\R$, $a<b$, such that $g_k(x,u)=0$ for all $u$ outside the compact $[a,b]$, for all $x\in\T^N$, $k\geq 1$. For simplicity, we will take $a=-1$, $b=1$. We will assume therefore that
\begin{equation}\label{multiplicativeNoise}
\mbox{for all }u\in\R,\;|u|\geq 1\Rightarrow g_k(x,u)=0,
\end{equation}
for all $x\in\T^N$, $k\geq 1$, and consider initial data with values in $[-1,1]$. The solution of the continuous equation \eqref{stoSCL} then takes values in $[-1,1]$ almost-surely (see \cite[Theorem~22]{DottiVovelle16a}). There is no loss in generality in considering that $A$ is globally Lipschitz continuous then:
\begin{equation}\label{ALip}
\mathrm{Lip}(A):=\sup_{\xi\in\R}|A'(\xi)|<+\infty.
\end{equation} 
In that framework, we will build a stable and convergent approximation to \eqref{stoSCL} by an explicit-in-time Finite Volume method. 
Under \eqref{multiplicativeNoise}, it is also natural to assume 
\begin{equation}\label{D0plus}
\GG^2(x,u)=\|\Phi(x,u)\|_H^2=\sum_{k\geq 1}|g_k(x,u)|^2\leq D_0,
\end{equation}
which is of course stronger than \eqref{D0}. We may also perform the analysis of convergence of the Finite Volume method under \eqref{D0} instead of \eqref{D0plus}, but this puts exponential factors in various estimates, whereas these factors are close to $1$ in the real implementation of the scheme.
\bigskip

\textbf{Numerical approximation.} Let us give a brief summary of the theory of \eqref{stoSCL} and of its approximation.
Different approximation schemes to stochastically forced first-order conservation laws have already been analysed: time-discrete schemes, \cite{HoldenRisebro91,Bauzet15,KarlsenStorrosten2017}, space-discrete scheme \cite{KoleyMajeeVallet2016}, space-time Finite Volume discrete schemes: 
\begin{itemize}
\item in space dimension $1$, with strongly monotone fluxes, \cite{KrokerRohde12}
\item in space dimension $N\geq 1$, by a flux-splitting scheme, \cite{BauzetCharrierGallouet14a},
\item in space dimension $N\geq 1$, for general schemes with monotone fluxes, \cite{BauzetCharrierGallouet14b},
\end{itemize}
The Cauchy or the Cauchy-Dirichlet problem associated to \eqref{stoSCL} have been studied in \cite{EKMS00,Kim03,FengNualart08,ValletWittbold09,DebusscheVovelle10,ChenDingKarlsen12,BauzetValletWittbold12,
BauzetValletWittbold14,KarlsenStorrosten15}.
\medskip

The approximation of scalar conservation laws with stochastic flux has also been considered in \cite{GessPerthameSouganidis2016} (time-discrete scheme) and
\cite{MohamedSeaidZahri13} (space discrete scheme). For the corresponding Cauchy Problem, see \cite{LionsPerthameSouganidis13a,LionsPerthameSouganidis13,LionsPerthameSouganidis14,
GessSouganidis2015,GessSouganidis2014,Hofmanova2016}.
\medskip

\textbf{Kinetic formulation.} To prove the convergence of the Finite Volume method with monotone fluxes, we will use the companion paper \cite{DottiVovelle16a} and a kinetic formulation of the Finite Volume scheme. The subject of \cite{DottiVovelle16a} is the convergence of approximations to \eqref{stoSCL} in the context of the kinetic formulation of scalar conservation laws. Such kinetic formulations have been developed in \cite{LionsPerthameTadmor94,LionsPerthameTadmor94ki,MakridakisPerthame03,Perthame98,PerthameBook}. In \cite{MakridakisPerthame03}, a kinetic formulation of Finite Volume E-schemes is given (and applied in particular to obtain sharp CFL criteria). For Finite Volume schemes with monotone fluxes, the kinetic formulation is simpler, we give it explicitly in Proposition~\ref{prop:MP03}. Based on the kinetic formulation and an energy estimate, we derive some a priori bounds on the numerical approximation (theses are ``weak $BV$ estimates" in the terminology of \cite[Lemma 25.2]{EGH00}), see Section~\ref{secEnergyEstimates}. These estimates are used in the proof of consistency of the scheme when we show that it gives rise to an approximate solution to \eqref{stoSCL} in the sense of Definition~\ref{def:appsol}. Our final result, \textit{cf.} Theorem~\ref{th:mainthm}, should be compared to \cite[Theorem~2]{BauzetCharrierGallouet14b}. This latter gives the convergence of the Finite Volume method with monotone fluxes in a very similar context, under the slightly stronger hypothesis that the ratio of the time step $\Delta t$ with the spatial characteristic size $h$ of the mesh tends to $0$ when $h$ tends to $0$. \medskip

\textbf{Plan of the paper.} The plan of the paper is the following one. In the preliminary section~\ref{sec:solutions}, we give a brief summary of the notion of solution and approximate solution to \eqref{stoSCL} developed in \cite{DottiVovelle16a}. In Section~\ref{secFVscheme} is described the kind of approximation to \eqref{stoSCL} by the Finite Volume method which we consider here. 
In Section~\ref{secKiFVscheme} we establish the kinetic formulation of the scheme. This numerical kinetic formulation is analysed as follows: energy estimates are derived in Section~\ref{secEnergyEstimates}, then we show in Section~\ref{sec:appki} that this gives rise to an approximate generalized solution in the sense of Definition~\ref{def:appsol}. We show some additional estimates and then conclude to the convergence of the scheme in Section~\ref{sec:Convergence}. This result is stated in Theorem~\ref{th:mainthm}.

\section{Generalized solutions, approximate solutions}\label{sec:solutions}

The object of this section is to recall several results concerning the solutions to the Cauchy Problem associated to \eqref{stoSCL} and their approximations. We give the main statements, without much explanations or comments; those latter can be found in \cite{DottiVovelle16a}: we give the precise references when needed.

\subsection{Solutions}\label{sec:solutionssolutions}

\begin{definition}[Random measure] Let $X$ be a topological space. If $m$ is a map from $\Omega$ into the set of non-negative finite Borel measures on $X$ such that, for each continuous and bounded function $\phi$ on $X$, $\<m,\phi\>$ is a random variable, then we say that $m$ is a random measure on $X$.
\label{def:randommeasure}\end{definition}

A random measure $m$ is said to have a finite first moment if
\begin{equation}\label{Firstm}
\E m(\T^N\times[0,T]\times\R)<+\infty.
\end{equation} 

\begin{definition}[Solution] Let $u_{0}\in L^\infty(\T^N)$. An $L^1(\T^N)$-valued stochastic process $(u(t))_{t\in[0,T]}$ is said to be a solution to~\eqref{stoSCL} with initial datum $u_0$ if $u$ and $\mathtt{f}:=\mathbf{1}_{u>\xi}$ have the following properties:
\begin{enumerate}
\item\label{item:1defsol} $u\in L^1_{\mathcal{P}}(\T^N\times [0,T]\times\Omega)$,
\item\label{item:1bisdefsol} for all $\varphi\in C^1_c(\T^N\times\R)$, almost-surely, $t\mapsto\<\mathtt{f}(t),\varphi\>$ is c{\`a}dl{\`a}g,
\item\label{item:2defsol} for all $p\in[1,+\infty)$, there exists $C_p\geq 0$ such that
\begin{equation}
\E\left(\sup_{t\in[0,T]}\|u(t)\|_{L^p(\T^N)}^p\right)\leq C_p,
\label{eq:integrabilityu}\end{equation} 
\item\label{item:3defsol} there exists a random measure $m$ with first moment \eqref{Firstm}, such that for all $\varphi\in C^1_c(\T^N\times\R)$, for all $t\in[0,T]$,
\begin{multline}
\<\mathtt{f}(t),\varphi\>=\<\mathtt{f}_0,\varphi\>
+\int_0^t \<\mathtt{f}(s),a(\xi)\cdot\nabla\varphi\>ds\\
+\sum_{k\geq 1}\int_0^t\int_{\T^N}g_k(x,u(x,s))\varphi(x,u(x,s)) dx d\beta_k(s)\\
+\frac{1}{2}\int_0^t\int_{\T^N} \partial_\xi\varphi(x,u(x,s))\GG^2(x,u(x,s)) dx ds-m(\partial_\xi\varphi)([0,t]),
\label{eq:kineticupre}\end{multline}
a.s., where $\mathtt{f}_0(x,\xi)=\mathbf{1}_{u_0(x)>\xi}$, $\GG^2:=\sum_{k=1}^\infty |g_k|^2$ and $a(\xi):=A'(\xi)$.
\end{enumerate}
\label{defkineticsol}\end{definition}

In item~\ref{item:1defsol}, the index $\mathcal{P}$ in $u\in L^1_{\mathcal{P}}(\T^N\times [0,T]\times\Omega)$ means that $u$ is predictable. See \cite[Section~2.1.1]{DottiVovelle16a}. The function denoted $\mathtt{f}:=\mathbf{1}_{u>\xi}$ is given more precisely by
$$
(x,t,\xi)\mapsto\mathbf{1}_{u(x,t)>\xi}.
$$
This is the characteristic function of the subgraph of $u$. To study the stability of solutions, or the convergence of approximate solutions (there are two similar problems), we have to consider the stability of this property, the fact of being the `` characteristic function of the subgraph of a function". If $(u_n)$, is a sequence of functions, say on a finite measure space $X$, $p\in(1,\infty)$ and $(u_n)$ is bounded in $L^p(X)$, then there is a subsequence still denoted $(u_n)$ which converges to a function $u$ in $L^p(X)$-weak. Up to a subsequence, the sequence of equilibrium functions 
$\mathtt{f}_n:=\mathbf{1}_{u_n>\xi}$ is converging to a function $f$ in $L^\infty(X\times\R)$-weak star. The limit $f$ is equal to $\mathtt{f}:=\mathbf{1}_{u>\xi}$ if, and only if, $(u_n)$ is converging strongly, see \cite[Lemma~2.6]{DottiVovelle16a}. When strong convergence remains a priori unknown, the limit $f$ still keeps some structural properties. This is a kinetic function in the sense of Definition~\ref{def:kifunction} below, \cite[Corollary~2.5]{DottiVovelle16a}. Our notion of generalized solution is based on this notion.

\subsection{Generalized solutions}\label{sec:ges}

\begin{definition}[Young measure] Let $(X,\mathcal{A},\lambda)$ be a finite measure space. Let $\mathcal{P}_1(\R)$ denote the set of probability measures on $\R$. We say that a map $\nu\colon X\to\mathcal{P}_1(\R)$ is a Young measure on $X$ if, for all $\phi\in C_b(\R)$, the map $z\mapsto \nu_z(\phi)$ from $X$ to $\R$ is measurable. We say that a Young measure $\nu$ vanishes at infinity if, for every $p\geq 1$, 
\begin{equation}
\int_X\int_\R |\xi|^p d\nu_z(\xi)d\lambda(z)<+\infty.
\label{nuvanish}\end{equation}
\label{defYoung}\end{definition}

\begin{definition}[Kinetic function] Let $(X,\mathcal{A},\lambda)$ be a finite measure space. A measurable function $f\colon X\times\R\to[0,1]$ is said to be a kinetic function if there exists a Young measure $\nu$ on $X$ that vanishes at infinity such that, for $\lambda$-a.e. $z\in X$, for all $\xi\in\R$,
\begin{equation*}
f(z,\xi)=\nu_{z}(\xi,+\infty).
\end{equation*}
We say that $f$ is an {\rm equilibrium} if there exists a measurable function $u\colon X\to\R$ such that $f(z,\xi)=\mathtt{f}(z,\xi)=\mathbf{1}_{u(z)>\xi}$ a.e., or, equivalently, $\nu_z=\delta_{\xi=u(z)}$ for a.e. $z\in X$.
\label{def:kifunction}\end{definition}

\begin{definition}[Generalized solution]
\label{d4} Let $f_0\colon\T^N\times\R\to[0,1]$ be a kinetic function. An $L^\infty(\T^N\times\R;[0,1])$-valued process $(f(t))_{t\in[0,T]}$ is said to be a generalized solution to~\eqref{stoSCL} with initial datum $f_0$ if $f(t)$ and $\nu_t:=-\partial_\xi f(t)$ have the following properties:
\begin{enumerate}
\item\label{item:1d4} for all $t\in[0,T]$, almost-surely, $f(t)$ is a kinetic function, and, for all $R>0$, $f\in L^1_\mathcal{P}(\T^N\times(0,T)\times(-R,R)\times\Omega)$,
\item\label{item:1terd4} for all $\varphi\in C^1_c(\T^N\times\R)$, almost-surely, the map $t\mapsto\<f(t),\varphi\>$ is c{\`a}dl{\`a}g,
\item\label{item:2d4} for all $p\in[1,+\infty)$, there exists $C_p\geq 0$ such that 
\begin{equation}
\E\left(\sup_{t\in[0,T]}\int_{\T^N}\int_\R|\xi|^p d\nu_{x,t}(\xi) dx\right) \leq C_p,
\label{eq:integrabilityf}\end{equation}
\item\label{item:4d4} there exists a random measure $m$  with first moment \eqref{Firstm}, such that for all $\varphi\in C^1_c(\T^N\times\R)$, for all $t\in[0,T]$, almost-surely,
\begin{align}
\<f(t),\varphi\>=&\<f_0,\varphi\>+\int_0^t \<f(s),a(\xi)\cdot\nabla_x\varphi\>ds
\nonumber\\
&+\int_0^t\int_{\T^N}\int_\R g_k(x,\xi)\varphi(x,\xi)d\nu_{x,s}(\xi) dxd\beta_k(s) \nonumber\\
&+\frac{1}{2}\int_0^t\int_{\T^N}\int_\R \GG^2(x,\xi)\partial_\xi\varphi(x,\xi)d\nu_{x,s}(\xi) dx ds
-m(\partial_\xi\varphi)([0,t]).
\label{eq:kineticfpre}
\end{align}
\end{enumerate}
\end{definition}

The following statement is Theorem~3.2. in \cite{DottiVovelle16a}.

\begin{theorem}[Uniqueness, Reduction] Let $u_0\in L^\infty(\T^N)$. Assume~\eqref{D0}-\eqref{D1}. 
Then we have the following results:
\begin{enumerate}
\item there is at most one solution with initial datum $u_0$ to \eqref{stoSCL}. 
\item If $f$ is a ge\-ne\-ra\-lized solution to \eqref{stoSCL} with initial 
datum $f_0$ \emph{at equilibrium:} $f_0=\mathbf{1}_{u_0>\xi}$, then there exists a solution $u$ to \eqref{stoSCL} with initial datum 
$u_0$ such that $f(x,t,\xi)=\mathbf{1}_{u(x,t)>\xi}$ a.s., for a.e. $(x,t,\xi)$. 
\item if $u_1$, $u_2$ are two solutions to \eqref{stoSCL} associated to the initial data $u_{1,0}$, $u_{2,0}\in L^\infty(\T^N)$ respectively, then
\begin{equation}
\E\|(u_1(t)-u_2(t))^+\|_{L^1(\T^N)}\leq\E\|(u_{1,0}-u_{2,0})^+\|_{L^1(\T^N)}.
\label{L1compadd}\end{equation}
This implies the $L^1$-contraction property, and comparison principle for solutions.
\end{enumerate} 
\label{th:Uadd}\end{theorem}

\subsection{Approximate solutions}\label{sec:solutionsAppsolutions}
In \cite{DottiVovelle16a}, we prove the convergence of solutions to different problems, which are approximations of \eqref{stoSCL}. The solutions of these approximate problems give rise to approximate solutions and, more precisely, of approximate generalized solutions, according to the following definition (see Definition~4.1 and Section~5 in \cite{DottiVovelle16a}).
\begin{definition}[Approximate generalized solutions] Let $f_0^n\colon\T^N\times\R\to[0,1]$ be some kinetic functions. Let $(f^n(t))_{t\in[0,T]}$ be a sequence of $L^\infty(\T^N\times\R;[0,1])$-valued processes. Assume that the functions $f^n(t)$, and the associated Young measures $\nu^n_t=-\partial_\xi\varphi f^n(t)$ are satisfying 
item \ref{item:1d4}, \ref{item:1terd4}, \ref{item:2d4}, in Definition~\ref{d4} and Equation~\eqref{eq:kineticfpre} up to an error term, \textit{i.e.}: for all $\varphi\in C^\order_c(\T^N\times\R)$, there exists an adapted process $\eps^n(t,\varphi)$, with $t\mapsto\eps^n(t,\varphi)$ almost-surely continuous such that
\begin{equation}\label{epsto0}
\lim\limits_{n\to+\infty}\sup_{t\in[0,T]}\left|\eps^n(t,\varphi)\right|=0\mbox{ in probability,}
\end{equation} 
and there exists some random measures $m^n$  with first moment \eqref{Firstm}, such that, for all $n$, for all $\varphi\in C^\order_c(\T^N\times\R)$, for all $t\in[0,T]$, almost-surely,
\begin{align}
\<f^n(t),\varphi\>=&\eps^n(t,\varphi)+\<f^n_0,\varphi\>+\int_0^t \<f^n(s),a(\xi)\cdot\nabla_x\varphi\>ds
\nonumber\\
&+\int_0^t\int_{\T^N}\int_\R g_k(x,\xi)\varphi(x,\xi)d\nu^n_{x,s}(\xi) dxd\beta_k(s) \nonumber\\
&+\frac{1}{2}\int_0^t\int_{\T^N}\int_\R \GG^2(x,\xi)\partial_\xi\varphi(x,\xi)d\nu^n_{x,s}(\xi) dx ds
-m^n(\partial_\xi\varphi)([0,t]).
\label{eq:kineticfpreappt}
\end{align}
Then we say that $(f^n)$ is a sequence of approximate generalized solutions to~\eqref{stoSCL} with initial datum $f_0^n$.
\label{def:appsol}\end{definition}

Consider a sequence $(f_n)$ of approximate solutions to \eqref{stoSCL} satisfying the following (minimal) bounds.
\begin{enumerate}
\item There exists $C_p\geq 0$ independent on $n$ such that $\nu^n:=-\partial_\xi f^n$ satisfies 
\begin{equation}
\E\left[\sup_{t\in[0,T]}\int_{\T^N}\int_\R|\xi|^p d\nu^n_{x,t}(\xi) dx\right]\leq C_p,
\label{eq:integrabilityfn}\end{equation}
\item the measures $\E m^n$ satisfy the bound
\begin{equation}\label{Boundmn}
\sup_n\E m^n(\T^N\times[0,T]\times\R)<+\infty,
\end{equation}
and the following tightness condition: if $B_R^c=\{\xi\in\R,|\xi|\geq R\}$, then
\begin{equation}
\lim_{R\to+\infty}\sup_n\E m^n(\T^N\times[0,T]\times B_R^c)=0.
\label{inftymn}\end{equation}
\end{enumerate}

We give in \cite{DottiVovelle16a} the proof of the following convergence result, see Theorem~40 in \cite{DottiVovelle16a}.

\begin{theorem}[Convergence] Suppose that there exists a sequence of approximate ge\-ne\-ra\-li\-zed solutions $(f^n)$ to~\eqref{stoSCL} with initial datum $f_0^n$ satisfying \eqref{eq:integrabilityfn}, \eqref{Boundmn} and the tightness condition \eqref{inftymn} and such that $(f^n_0)$ converges to the equilibrium function $\mathtt{f}_0(\xi)=\mathbf{1}_{u_0>\xi}$ in $L^\infty(\T^N\times\R)$-weak-*, where $u_0\in L^\infty(\T^N)$. We have then
\begin{enumerate}
\item there exists a unique solution $u\in L^1(\T^N\times [0,T]\times\Omega)$ to \eqref{stoSCL} with initial datum $u_0$;
\item let
$$
u^n(x,t)=\int_\R\xi d\nu^n_{x,t}(\xi)=\int_\R\left( f^n(x,t,\xi)-\mathbf{1}_{0>\xi} \right)d\xi.
$$
Then, for all $p\in[1,\infty[$, $(u^n)$ is converging to $u$ with the following two different modes of convergence: $u_n\to u$ in $L^p(\T^N\times(0,T)\times\Omega)$ and almost surely, for all $t\in[0,T]$, $u_n(t)\to u(t)$ in $L^p(\T^N)$.
\end{enumerate}
\label{th:pathcv} \end{theorem}

In the next section, we define the numerical approximation to \eqref{stoSCL} by the Finite Volume method. To prove the convergence of the method, we will show that the hypotheses of Theorem~\ref{th:pathcv} are satisfied. The most difficult part in this programme is to prove that the numerical approximations generates a sequence of approximate generalized solutions, see Section~\ref{sec:appki}.

\section{The finite volume scheme}\label{secFVscheme}

\paragraph{Mesh} A mesh of $\T^N$ is a family ${\mathcal{T}/\Z^N}$ of disjoint connected open subsets $K\in(0,1)^N$ which form a partition of $(0,1)^N$ up to a negligible set. We denote by $\mathcal{T}$ the mesh
$$
\{K+l;l\in\Z^N,K\in{\mathcal{T}/\Z^N}\}
$$
deduced on $\R^N$. For all distinct $K,L\in\mathcal{T}$, we assume that $\overline{K}\cap\overline{L}$ is contained in an hyperplane; the interface between $K$ and $L$ is denoted $K|L:=\overline{K}\cap\overline{L}$. The set of neighbours of $K$ is
\begin{equation*}
\mathcal{N}(K)=\left\{L\in\mathcal{T}; L\not= K,\; K|L\not=\emptyset\right\}.
\end{equation*}
We use also the notation
\begin{equation*}
\partial K=\bigcup_{L\in\mathcal{N}(K)} K|L.
\end{equation*}
In general, there should be no confusion between $\partial K$ and the topological boundary $\overline{K}\setminus K$. 

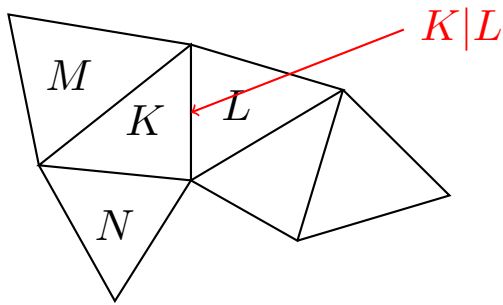
\begin{center}
\begin{tikzpicture}[thick,scale=2] 
\draw (0,0) --  (1,-0.1) -- (1,0.8) -- cycle;
\node[scale=1.5] at (0.7,0.3) {$K$}; 
\draw (1,-0.1) -- (2,0.5) -- (1,0.8);
\node[scale=1.5] at (1.3,0.4) {$L$};
\draw (0,0) -- (-0.2,1) -- (1,0.8);
\node[scale=1.5] at (0.2,0.6) {$M$};
\draw (0,0) -- (0.5,-0.9) -- (1,-0.1);
\node[scale=1.5] at (0.5,-0.4) {$N$}; 
\draw (1,-0.1) -- (1.7,-0.5) -- (2,0.5);
\draw (1.7,-0.5) -- (2.7,-0.2) -- (2,0.5);
\draw[red,->] (2.4,0.9) node[anchor=west,scale=1.5]{$K|L$} -- (1,0.35);
\end{tikzpicture}
\end{center}
\medskip

We also denote by $|K|$ the $N$-dimensional Lebesgue Measure of $K$ and by $|\partial K|$ (respectively $|K|L|$) the $(N-1)$-dimensional Haussdorff measure of $\partial K$ (respectively of $K|L$) (the $(N-1)$-dimensional Haussdorff measure is normalized to coincide with the $(N-1)$-dimensional Lebesgue measure on hyperplanes).

\paragraph{Scheme} Let $(A_{K\to L})_{K\in\mathcal{T},L\in\mathcal{N}(K)}$ be a family of monotone, Lipschitz continuous numerical flux, consistent with $A$: we assume that each function $A_{K\to L}$ satisfies the following properties.
\begin{itemize}
\item Monotony: $A_{K\to L}(v,w)\leq A_{K\to L}(v',w)$ for all $v$, $v'$, $w\in\R$ with $v\leq v'$ and $A_{K\to L}(v,w)\geq A_{K\to L}(v,w')$ for all $v$, $w$, $w'\in\R$ with $w\leq w'$.
\item Lipschitz regularity: there exists $L_A\geq 0$ such that
\begin{equation}\label{AALip}
|A_{K\to L}(v,w)-A_{K\to L}(v',w')|\leq |K|L| L_A,
\end{equation}
for all $v$, $v'$, $w$, $w'\in\R$.
\item Consistency: 
\begin{equation}\label{consistency}
A_{K\to L}(v,v)=\int_{K|L} A(v)\cdot n_{K,L} d\mathcal{H}^{N-1}=|K|L|A(v)\cdot n_{K,L},
\end{equation}
for all $v\in\R$, where $n_{K,L}$ is the outward unit normal to $K$ on $K|L$.
\item Conservative symmetry:
\begin{equation}\label{conservativesym}
A_{K\to L}(v,w)=-A_{L\to K}(w,v),
\end{equation}
for all $K,L\in\mathcal{T}$, $v,w\in\R$.
\end{itemize}

The conservative symmetry property ensures that the numerical flux $Q^n_{K\to L}$ defined below in \eqref{Halfflux} satisfies $Q^n_{K\to L}=-Q^n_{L\to K}$ for all $K,L$. \medskip

Let $t_n<t_{n+1}$ be two given discrete time. Let $\Delta t_n=t_{n+1}-t_n$. Knowing $v^n_K$, an approximation of the value of the solution $u$ to \eqref{stoSCL} in the cell $K$ at time $t_n$, we compute $v^{n+1}_K$, the approximation to the value of $u$ in $K$ at the next time step $t_{n+1}$, by the formula
\begin{equation}\label{FVscheme}
|K|(v^{n+1}_K-v^n_K)+\Delta t_n\sum_{L\in\mathcal{N}(K)}Q^n_{K\to L}=|K| (\Delta t_n)^{1/2} g_{k,K}(v^n_K)X^{n+1}_k,
\end{equation}
where $K\in\mathcal{T}$, with the initialization
\begin{equation}\label{FVIC}
v^0_K=\frac{1}{|K|}\int_{K}u_0(x)dx,\quad K\in\mathcal{T}.
\end{equation} 
In \eqref{FVscheme}, $\Delta t_n Q^n_{K\to L}$ is the numerical flux at the interface $K|L$ on the range of time $[t_n,t_{n+1}]$, where $Q^n_{K\to L}$ is given by
\begin{equation}\label{Halfflux}
Q^n_{K\to L}=A_{K\to L}(v^n_K,v^n_L).
\end{equation}
We have also defined
\begin{equation}\label{defX}
X^{n+1}_k=\frac{\beta_k(t_{n+1})-\beta_k(t_n)}{(\Delta t_n)^{1/2}}.
\end{equation}
Then, the $(X^{n+1}_k)_{k\geq 1,n\in\N}$ are i.i.d. random variables with normalized centred normal law $\mathcal{N}(0,1)$. Besides, for each $n\geq 1$, the sequence $(X^{n+1}_k)_{k\geq 1}$ is independent on $\mathcal{F}_n$, the sigma-algebra generated by $\{X^{m+1}_k;k\geq 1,m<n\}$. The numerical functions $g_{k,K}$ are defined by the average
\begin{equation}\label{defgnum}
g_{k,K}(v)=\frac{1}{|K|}\int_{K} g_k(x,v) dx.
\end{equation}
Then, in virtue of \eqref{D0plus} we have
\begin{equation}
\GG^2_K(v):=\sum_{k\geq 1}|g_{k,K}(v)|^2\leq D_0, \label{D0num}
\end{equation}
where $v\in\R$, $K\in\mathcal{T}$. We deduce \eqref{D0num} from \eqref{D0plus} and Jensen's Inequality. Similarly, we deduce from \eqref{D1} and Jensen's Inequality that
\begin{equation*}
\sum_{k\geq 1}|g_{k,K}(\xi)-g_k(y,\xi)|^2\leq D_1\frac{1}{|K|}\int_K |x-y|^2 dx,
\end{equation*}
for all $y\in\T^N$. In particular (switching from the variable $y$ to the variable $x$), and assuming $\mathrm{diam}(K)\leq h$ (this is the hypothesis \eqref{hsizemesh}, which we will make later), we have the following consistency estimate
\begin{equation}\label{D1num0}
\sum_{k\geq 1}|g_{k,K}(\xi)-g_k(x,\xi)|^2\leq D_1 h^2,
\end{equation}
for all $x\in K$, which will be used later (see \eqref{D0D1} for example).

\begin{remark}[Approximation in law] In effective computations, the random variables $X^{n+1}_k$ are drawn at each time step. They are i.i.d. random variables with normalized centred normal law $\mathcal{N}(0,1)$. In this situation, we will prove the convergence in law of the Finite Volume scheme to the solution to \eqref{stoSCL}, see Remark~\ref{rk:rkcvlaw} after Theorem~\ref{th:mainthm}.
\end{remark}

\begin{remark}[Global Lipschitz Numerical Flux] We assume in \eqref{AALip} that the numerical fluxes $A_{K\to L}$ are globally Lipschitz continuous. This is consistent with \eqref{ALip}. Both \eqref{AALip} and \eqref{ALip} are strong hypotheses, except if a priori $L^\infty$-bounds are known on the solutions to \eqref{stoSCL}, which is the case here, thanks to the hypothesis of compact support \eqref{multiplicativeNoise}. Without loss of generality, we will assume that $\mathrm{Lip}(A)\leq L_A$.
\label{rk:LA}\end{remark}
\section{The kinetic formulation of the finite volume scheme}\label{secKiFVscheme}

The kinetic formulation of the Finite Volume method has been introduced by Makridakis and Perthame in \cite{MakridakisPerthame03}. The principle is the following one. For linear transport equations, which corresponds to a linear flux function $A(u)=au$, $a\in\R^N$, the upwind numerical flux $A_{K\to L}$ in \eqref{Halfflux} is given by 
\begin{equation}\label{EOtransport}
A_{K\to L}(v,w)=[a^*_{K\to L}]^+ v-[a^*_{K\to L}]^- w,
\end{equation}
where
\begin{equation}\label{astarKL}
a^*_{K\to L}:=\int_{K|L}a\cdot n_{K|L}d\mathcal{H}^{N-1}=|K|L|a\cdot n_{K,L},
\end{equation}
with the usual notation $v^+=\max(v,0)$, $w^-=(-w)^+$. The discrete approximation of the transport equation
$$
\partial_t f+a(\xi)\cdot\nabla_x f=0
$$
by the Finite Volume method is therefore
\begin{equation}\label{upwindKi}
|K|(f^{n+1}_K-f^n_K)+\Delta t_n \sum_{L\in\mathcal{N}(K)} a_{K\to L}^n(\xi)=0,
\end{equation}
where
$$
a_{K\to L}^n(\xi)=[a^*_{K|L}(\xi)]^+f^n_K-[a^*_{K|L}(\xi)]^-f^n_{L}.
$$
Recall that $\GG_K$ is defined by \eqref{D0num}. A kinetic formulation of \eqref{FVscheme} consistent with \eqref{upwindKi} would be
\begin{multline}\label{KiFVNo}
|K|(\mathtt{f}^{n+1}_K-\mathtt{f}^n_K)+\Delta t_n \sum_{L\in\mathcal{N}(K)} a_{K\to L}^n(\xi)\\
=|K| (\Delta t_n)^{1/2}\delta_{v^n_K=\xi} g_{k,K}(\xi)X^{n+1}_k+|K|\Delta t_n\partial_\xi\Big(m^n_K(\xi)-\frac12\GG^2_K(\xi)\delta_{v^n_K=\xi}\Big),
\end{multline}
where, for $K\in\mathcal{T}$, $n\in\N$, $\xi\in\R$, $m^n_K(\xi)\geq 0$, and where
\begin{equation}\label{equilibriumFV}
\mathtt{f}^n_K(\xi):=\mathbf{1}_{v^n_K>\xi}.
\end{equation}
This is not exactly the kinetic formulation that we will consider. See \eqref{eq:discretekieq} for a correct version of \eqref{KiFVNo}. We will mainly work with a kinetic formulation (see \eqref{KiFV}), obtained thanks to the following splitting method. For $K\in\mathcal{T}$ and $n\in\N$, let us define $v^{n+1/2}_K$ as the solution to   
\begin{equation}\label{vhalf}
|K|(v^{n+1/2}_K-v^n_K)+\Delta t_n\sum_{L\in\mathcal{N}(K)}A_{K\to L}(v^n_K,v^n_L)=0.
\end{equation}
Then $v^{n+1/2}_K$ is the state reached after a step of deterministic evolution, by the discrete approximation of the equation $u_t+\div(A(u))=0$. To this step corresponds the kinetic formulation
\begin{equation}\label{KiFVhalf}
|K|(\mathtt{f}^{n+1/2}_K-\mathtt{f}^n_K)+\Delta t_n \sum_{L\in\mathcal{N}(K)} a_{K\to L}^n(\xi)=|K|\Delta t_n\, \partial_\xi m^n_K(\xi),
\end{equation}
where $\mathtt{f}^m_K(\xi)=\mathbf{1}_{v^m_K>\xi}$, $m\in\{n,n+1/2\}$ and 
\begin{equation}\label{mnkpositive}
m^n_K\geq 0.
\end{equation} 
In \eqref{KiFVhalf}, $a_{K\to L}^n(\xi)$ is a function
\begin{equation}
a_{K\to L}^n(\xi)=a_{K\to L}(\xi,v_K^n,v_L^n),
\end{equation}
where $(\xi,v,w)\mapsto a_{K\to L}(\xi,v,w)$ satisfies the following consistency conditions:
\begin{align}
\int_\R \left[a_{K\to L}(\xi,v,w)-a^*_{K\to L}(\xi)\mathbf{1}_{0>\xi}\right] d\xi=A_{K\to L}(v,w), \label{consistencya1}\\
a_{K\to L}(\xi,v,v)=a^*_{K\to L}(\xi)\mathbf{1}_{v>\xi},\label{consistencya2}
\end{align}
for all $\xi$, $v$, $w\in\R$, where $a^*_{K\to L}$ is defined by \eqref{astarKL}. Before we prove the existence of the kinetic formulation \eqref{KiFVhalf}-\eqref{mnkpositive}-\eqref{consistencya1}-\eqref{consistencya2}, see Proposition~\ref{prop:MP03}, let us first deduce from \eqref{KiFVhalf} the kinetic formulation of the whole scheme \eqref{FVscheme}. This is the equation
\begin{multline}\label{KiFV}
|K|(\mathtt{f}^{n+1}_K(\xi)-\mathtt{f}^n_K(\xi))+\Delta t_n \sum_{L\in\mathcal{N}(K)} a_{K\to L}^n(\xi)\\
=|K|\Delta t_n\, \partial_\xi m^n_K(\xi)+|K|\Big[\mathtt{f}^{n+1}_K(\xi)-\mathtt{f}^{n+1/2}_K(\xi)\Big].
\end{multline}
We may try to develop the term $\mathtt{f}^{n+1}_K(\xi)-\mathtt{f}^{n+1/2}_K(\xi)$ (this is done in \eqref{discretekieq4} to obtain \eqref{eq:discretekieq}), but \eqref{KiFV} will be sufficient for the moment. It will be sufficient in particular to obtain the so-called energy estimates of Section~\ref{secEnergyEstimates}.

\begin{proposition}[Kinetic formulation of the Finite Volume method]\label{prop:MP03} Set 
\begin{equation}
a_{K\to L}(\xi,v,w)=a^*_{K\to L}(\xi)\mathbf{1}_{\xi<v\wedge w}+\left[\partial_2A_{K\to L}(v,\xi)\mathbf{1}_{v\leq \xi\leq w}
+\partial_1 A_{K\to L}(\xi,w)\mathbf{1}_{w\leq \xi\leq v}\right]
\label{newamon}
\end{equation}
and
\begin{equation}\label{defmnK}
m^n_K(\xi)=-\frac{1}{\Delta t_n}\left[(v^{n+1/2}_K-\xi)^+-(v^n_K-\xi)^+\right]-\frac{1}{|K|}\sum_{L\in\mathcal{N}(K)}\int_\xi^{+\infty}a^n_{K\to L}(\zeta)d\zeta.
\end{equation}
Let us also assume that 
\begin{equation}\label{CFLplus}
\Delta t_n\frac{|\partial K|}{|K|}L_A\leq 1,\quad\forall K\in\mathcal{T},
\end{equation}
for all $n\in\N$, $K\in\mathcal{T}$. Then the equations \eqref{KiFVhalf}-\eqref{mnkpositive}-\eqref{consistencya1}-\eqref{consistencya2} are satisfied and, besides, we have
\begin{equation}\label{supportaKLxi}
a_{K\to L}(\xi,v,w)=0,\quad\mbox{when }\xi\geq v\vee w,
\end{equation}
for all $K,L\in\mathcal{T}$.
\end{proposition}

\begin{remark}[Support of $m^n_K$] By \eqref{defmnK}, the definition \eqref{newamon} of $a_{K\to L}$ and the equation \eqref{vhalf}, $\xi\mapsto m^n_K(\xi)$ is compactly supported in the convex envelope of $v^{n+1/2}_K$, $v^n_K$, $\{v^n_L;L\in\mathcal{N}(K)\}$.
\label{rk:suppmnK}\end{remark}

\textbf{Proof of Proposition~\ref{prop:MP03}.} We check at once \eqref{KiFVhalf} and \eqref{consistencya1}, \eqref{consistencya2}, \eqref{supportaKLxi}. To show that $m^n_K(\xi)\geq 0$, let us introduce
\begin{align}
\Phi_{K\to L}(\xi,v,w)&=\int_\xi^{+\infty} a_{K\to L}(\zeta,v,w)d\zeta\label{defPhinum}\\
\Phi^n_{K\to L}(\xi)&=\Phi_{K\to L}(\xi,v^n_K,v^n_L)=\int_\xi^{+\infty} a^n_{K\to L}(\zeta)d\zeta.\label{defPhinumn}
\end{align}
A simple computation gives the formula
\begin{equation}\label{eq:numEntropyFlux}
\Phi^n_{K\to L}(\xi)=A_{K\to L}(v^n_K,v^n_L)-A_{K\to L}(v^n_K\wedge\xi,v^n_L\wedge\xi).
\end{equation}
By comparison with the identity $(v-\xi)^+=v-v\wedge\xi$, the quantities $\Phi^n_{K\to L}(\xi)$ appears, in virtue of \eqref{eq:numEntropyFlux}, as the numerical entropy fluxes associated to the entropy $\eta(v):=(v-\xi)^+$. Then $m^n_K(\xi)\geq 0$ is equivalent to the discrete entropy inequality
\begin{equation}\label{eq:discreteEntropyInequality}
\frac{1}{\Delta t_n}\left[\eta(v^{n+1/2}_K)-\eta(v^n_K)\right]+\frac{1}{|K|}\sum_{L\in\mathcal{N}(K)}\Phi^n_{K\to L}(\xi)\leq 0.
\end{equation}
It is a classical fact that, under the CFL condition \eqref{CFLplus}, the deterministic Finite Volume scheme~\eqref{vhalf} has the following monotony property: $v^{n+1/2}_K$ in \eqref{vhalf} is a non-decreasing function of each of the entries $v^n_K$, $v^n_L$, $L\in\mathcal{N}(K)$. This implies \eqref{eq:discreteEntropyInequality} then. See Lemma~25.1 and Lemma~27.1 in \cite{EGH00}. \qed

\section{Energy estimates}\label{secEnergyEstimates}

The Finite Volume scheme \eqref{FVscheme} may be compared to the stochastic parabolic equation
\begin{equation}\label{StoParabolic}
du^\eps(x,t)+\div(A(u^\eps(x,t)))dt-\eps\Delta u^\eps(x,t) dt=\Phi(x,u^\eps(x,t))dW(t).
\end{equation}
For \eqref{StoParabolic}, we have the energy estimate
\begin{equation}\label{EnergyEstimate}
\frac12\frac{d\;}{dt}\E\|u^\eps(t)\|_{L^2(\T^N)}^2+\eps\E\|\nabla u^\eps(t)\|_{L^2(\T^N)}^2=\frac12\E\|\GG(\cdot,u^\eps(\cdot,t))\|_{L^2(\T^N)}^2.
\end{equation}
(Recall that $\GG$ is defined by \eqref{D0}). In the following Proposition~\ref{prop:EnergyEstimate}, we obtain an analogous result for the Finite Volume scheme \eqref{FVscheme}. To state Proposition~\ref{prop:EnergyEstimate}, we need first some notations.

\subsection{Notations}\label{sec:EEnotations}

For a fixed final time $T>0$, we denote by $\mathfrak{d}_T$ the set of admissible space-step and time-steps, defined as follows: if $h>0$ and $(\Delta t)=(\Delta t_0,\ldots,\Delta t_{N_T-1})$, $N_T\in\N^*$, then we say that $\delta:=(h,(\Delta t))\in \mathfrak{d}_T$ if 
\begin{equation}\label{ADMdelta1}
\frac{1}{h}\in\N^*,\quad t_{N_T}:=\sum_{n=0}^{N_T-1}\Delta t_n=T,\quad\sup_{0\leq n<N_T}\Delta t_n\leq 1.
\end{equation}
We say that $\delta\to 0$ if
\begin{equation}\label{ADMdelta3}
|\delta|:=h+\sup_{0\leq n<N_T}\Delta t_n\to 0.
\end{equation}
For a given mesh parameter $\delta=(h,(\Delta t))\in \mathfrak{d}_T$, we assume that a mesh $\mathcal{T}$ is given, with the following properties:
\begin{align}
\mathrm{diam}(K)\leq h,\label{hsizemesh}\\
\alpha_N h^N\leq |K|,\label{alphaK}\\
\ds|\partial K|\leq\frac{1}{\alpha_N} h^{N-1},\label{alphapK}
\end{align}
for all $K\in\mathcal{T}$, where 
$$
\mathrm{diam}(K)=\max_{x,y\in K}|x-y|
$$
is the diameter of $K$ and $\alpha_N$ is a given positive absolute constant depending on the dimension $N$ only. Note the following consequence of \eqref{alphaK}-\eqref{alphapK}:
\begin{equation}\label{alphaKpK}
h|\partial K|\leq\frac{1}{\alpha_N^2}|K|,
\end{equation}
for all $K\in\mathcal{T}$. We introduce then the discrete unknown $v_\delta(t)$ defined a.e. by
\begin{equation}\label{defvh}
v_\delta(x,t)=v^n_K,\quad x\in K, t_n\leq t<t_{n+1}.
\end{equation}
We will also need the intermediary discrete function
\begin{equation}\label{vflat}
v_\delta^\flat(x,t_{n+1})=v^{n+1/2}_K,\quad x\in K,
\end{equation}
defined for $n\in\N$.
Let us define the conjugate function $\bar f=1-f$. We introduce the following conjugate quantities:
\begin{equation}\label{conjugatea}
\bar a_{K\to L}(\xi,v,w)=a^*_{K\to L}(\xi)-a_{K\to L}(\xi,v,w),\quad\overline{\Phi}_{K\to L}(\xi,v,w)=\int_{-\infty}^\xi \bar a_{K\to L} (\zeta,v,w) d\zeta.
\end{equation}
We compute
\begin{equation}\label{barA0}
\overline{\Phi}_{K\to L}(\xi,v,w)=A_{K\to L}(\xi,\xi)-A_{K\to L}(v\wedge\xi,w\wedge\xi).
\end{equation}
We recognize in \eqref{barA0} a numerical flux associated to the entropy
$$
v\mapsto(v-\xi)^-=\xi-v\wedge\xi.
$$
From the explicit formula \eqref{newamon}, we obtain the identity 
\begin{multline}
\bar a_{K\to L}(\xi,v,w)=a_{K\to L}^*(\xi)\mathbf{1}_{\xi>v\vee w}+(a_{K\to L}^*(\xi)-\partial_2A_{K\to L}(v,\xi))\mathbf{1}_{v\leq \xi\leq w}\\
+(a_{K\to L}^*(\xi)-\partial_1A_{K\to L}(\xi,w))\mathbf{1}_{w\leq \xi\leq v}.
\label{newbaramon}
\end{multline}
Note that, for $a_{K\to L}$ defined by \eqref{newamon}, we have (using the fact that $\mathrm{Lip}(A)\leq L_A$),
\begin{equation}\label{BOUNDa}
\sup\{|a_{K\to L}(\xi,v,w)|;\xi,v,w\in\R\}\leq {L_A}|K|L|.
\end{equation}
Formula \eqref{newbaramon} gives the estimate
\begin{equation}\label{BOUNDbara}
\sup\{|\bar a_{K\to L}(\xi,v,w)|;\xi,v,w\in\R\}\leq 2{L_A}|K|L|,
\end{equation}
which is not optimal as \eqref{BOUNDa} may be, since it has an additional factor $2$. Consequently, we will use a slightly different formulation for $\overline{\Phi}_{K\to L}$:
\begin{align}
\overline{\Phi}_{K\to L}(\xi,v,w)&=\int_{-\infty}^\xi \bar{b}_{K\to L} (\zeta,\xi,v,w) d\zeta,\label{newbarA0}
\end{align}
where
\begin{align}
\bar{b}_{K\to L}(\zeta,\xi,v,w)&:=a^*_{K\to L}(\xi)\mathbf{1}_{\xi> v\vee w}+\partial_1 A_{K\to L}(\zeta,\xi)\mathbf{1}_{v\leq \xi\leq w}+\partial_2 A_{K\to L}(\xi,\zeta)\mathbf{1}_{w\leq \xi\leq v}.\label{defastar}
\end{align}
We also introduce
\begin{equation}\label{defbarbnKL}
\bar{b}_{K\to L}^n(\zeta,\xi)=\bar{b}_{K\to L}(\zeta,\xi,v^n_K,v^n_L).
\end{equation}
Now for $\bar{b}_{K\to L}$, we have an estimate similar to \eqref{BOUNDa}:
\begin{equation}\label{BOUNDastar}
\sup\{|\bar{b}_{K\to L}(\xi,v,w)|;\xi,\zeta,v,w\in\R\}\leq {L_A}|K|L|.
\end{equation}

\subsection{Energy estimate and controls by the dissipation}\label{sec:EEestimate}

\begin{proposition}[Energy estimate for the Finite Volume Scheme] Let $u_0\in L^\infty(\T)$, $T>0$ and $\delta\in \mathfrak{d}_T$. Let $(v_\delta(t))$ be the numerical unknown defined by \eqref{FVscheme}-\eqref{FVIC}-\eqref{defvh}. Set
\begin{equation}\label{EDissipation}
\mathcal{E}(T)=\sum_{n=0}^{N_T-1}\Delta t_n\sum_{K\in{\mathcal{T}/\Z^N}}|K|\int_\R m^n_K(\xi)d\xi.
\end{equation}
Then, under the CFL condition \eqref{CFLplus}, we have the energy estimate
\begin{equation}\label{eq:FVEnergyEestimate}
\frac12\E\|v_\delta(T)\|_{L^2(\T^N)}^2+\E\mathcal{E}(T)=\frac12\|v_\delta(0)\|_{L^2(\T^N)}^2+\frac12\E \sum_{n=0}^{N_T-1}\Delta t_n\sum_{K\in{\mathcal{T}/\Z^N}}|K|\sum_{k\geq 1} |g_{k,K}(v^n_K)|^2.
\end{equation}
\label{prop:EnergyEstimate}\end{proposition}

In the following proposition we derive various estimates, where the right-hand side is controlled by the dissipation term $\mathcal{E}(T)$ introduced in \eqref{EDissipation}.

\begin{proposition}[Control by the dissipation] Let $u_0\in L^\infty(\T^N)$, $T>0$ and $\delta\in \mathfrak{d}_T$. Let $(v_\delta(t))$ be the numerical unknown defined by \eqref{FVscheme}-\eqref{FVIC}-\eqref{defvh}. Let $v_\delta^\flat$ be defined by \eqref{vflat}. Then, under the CFL condition
\begin{equation}\label{CFLstrong}
2\Delta t_n\frac{|\partial K|}{|K|}\sup_{\xi\in\R}\frac{|a_{K\to L}^n(\xi)|}{|K|L|}\leq (1-\theta),\quad 0\leq n<N,\; K,L\in\mathcal{T},
\end{equation}
where $\theta\in(0,1)$, we have the following control:
\begin{align}
\sum_{n=0}^{N_T-1}\Delta t_n\sum_{K\in{\mathcal{T}/\Z^N}}\sum_{L\in\mathcal{N}(K)} &
\int_\R (\bar{\mathtt{f}}^n_{L}-\bar{\mathtt{f}}^n_K)\Phi_{K\to L}^n(\xi) d\xi\nonumber\\
\leq & 
\frac{2}{\theta}\,\sum_{n=0}^{N_T-1}\Delta t_n\sum_{K\in{\mathcal{T}/\Z^N}}|K| \int_\R \bar{\mathtt{f}}^n_K(\xi) m^n_K(\xi) d\xi,\label{eq:ControlSpace}
\end{align}
and
\begin{equation}\label{eq:ControlTime}
\sum_{n=0}^{N_T-1}\left\|\big[v_\delta^\flat(t_{n+1})-v_\delta(t_n)\big]_+\right\|_{L^2(\T^N)}^2
\leq 
\frac{2}{\theta}\,\sum_{n=0}^{N_T-1}\Delta t_n\sum_{K\in{\mathcal{T}/\Z^N}}|K| \int_\R \bar{\mathtt{f}}^n_K(\xi) m^n_K(\xi) d\xi.
\end{equation}
Under the CFL condition
\begin{equation}\label{CFLstrongbar}
2\Delta t_n\frac{|\partial K|}{|K|}\sup_{\xi\in\R}\frac{|\bar{b}_{K\to L}^n(\xi,\xi)|}{|K|L|}\leq (1-\theta),\quad 0\leq n<N,\; K,L\in\mathcal{T},
\end{equation}
where $\theta\in(0,1)$, (and where $\bar{b}_{K\to L}^n$ is defined by \eqref{defbarbnKL}) we have the following control:
\begin{align}
\sum_{n=0}^{N_T-1}\Delta t_n\sum_{K\in{\mathcal{T}/\Z^N}}\sum_{L\in\mathcal{N}(K)}  &
\int_\R (\mathtt{f}^n_{L}-\mathtt{f}^n_K)\overline{\Phi}_{K\to L}^n(\xi) d\xi\nonumber\\
\leq &
\frac{2}{\theta}\,\sum_{n=0}^{N_T-1}\Delta t_n\sum_{K\in{\mathcal{T}/\Z^N}}|K| \int_\R \mathtt{f}^n_K(\xi) m^n_K(\xi) d\xi,\label{eq:ControlSpacebar}
\end{align}
and
\begin{equation}\label{eq:ControlTimebar}
\sum_{n=0}^{N_T-1}\left\|\big[v_\delta^\flat(t_{n+1})-v_\delta(t_n)\big]_-\right\|_{L^2(\T^N)}^2
\leq 
\frac{2}{\theta}\,\sum_{n=0}^{N_T-1}\Delta t_n\sum_{K\in{\mathcal{T}/\Z^N}}|K| \int_\R \mathtt{f}^n_K(\xi) m^n_K(\xi) d\xi.
\end{equation}
\label{prop:ControlDissipation}\end{proposition}

Eventually, as a corollary to Proposition~\ref{prop:ControlDissipation}, we obtain the following estimates.
\begin{corollary}[Weak derivative estimates] Let $u_0\in L^\infty(\T^N)$, $T>0$ and $\delta\in \mathfrak{d}_T$. 
Assume that \eqref{multiplicativeNoise}, \eqref{D0plus}, \eqref{AALip}, \eqref{hsizemesh}, \eqref{alphaK} and \eqref{alphapK} are satisfied and that
\begin{equation}\label{CFLstrongU}
\Delta t_n \leq (1-\theta)\frac{\alpha_N^2}{2 {L_A}}\, h,\quad 0\leq n<N_T,
\end{equation}
where $\theta\in(0,1)$. Let $(v_\delta(t))$ be the numerical unknown defined by \eqref{FVscheme}-\eqref{FVIC}-\eqref{defvh}. Let $v_\delta^\flat$ be defined by \eqref{vflat}. Then we have the spatial estimate
\begin{multline}\label{eq:ControlSpaceAdd}
\E\sum_{n=0}^{N_T-1}\Delta t_n\sum_{K\in{\mathcal{T}/\Z^N}}\sum_{L\in\mathcal{N}(K)} \int_\R
\left[(\bar{\mathtt{f}}^n_{L}-\bar{\mathtt{f}}^n_K)\Phi_{K\to L}^n(\xi)
+ (\mathtt{f}^n_{L}-\mathtt{f}^n_K)\overline{\Phi}_{K\to L}^n(\xi)
\right]d\xi\\
\leq \frac{1}{\theta}\|v_\delta(0)\|_{L^2(\T^N)}^2+\frac{D_0 T}{\theta},
\end{multline}
and the two following temporal estimates:
\begin{equation}\label{eq:ControlTimeAdd}
\E\sum_{n=0}^{N_T-1}\left\|v_\delta^\flat(t_{n+1})-v_\delta(t_n)\right\|_{L^2(\T^N)}^2\leq 
\frac{1}{\theta}\|v_\delta(0)\|_{L^2(\T^N)}^2+\frac{D_0 T}{\theta},
\end{equation}
and
\begin{equation}\label{eq:ControlTimeAdd01}
\E\sum_{n=0}^{N_T-1}\left\|v_\delta(t_{n+1})-v_\delta(t_n)\right\|_{L^2(\T^N)}^2\leq 
\frac{1}{\theta}\|v_\delta(0)\|_{L^2(\T^N)}^2+\frac{2 D_0 T}{\theta}.
\end{equation}
\label{cor:EnergyEstimate}\end{corollary}

\subsection{Proof of Proposition~\ref{prop:EnergyEstimate}, Proposition~\ref{prop:ControlDissipation}, Corollary~\ref{cor:EnergyEstimate}}\label{sec:EEproofs}

\textbf{Proof of Proposition~\ref{prop:EnergyEstimate}.}  
We multiply first \eqref{KiFVhalf} by $\xi$ and sum the result over $K\in\mathcal{T}/\Z^N$ and $\xi\in\R$ to get the following balance equation 
\begin{equation}\label{eq:FVEnergyEestimate0}
\frac12\|v_\delta^\flat(t_{n+1})\|_{L^2(\T^N)}^2+\Delta t_n\sum_{K\in{\mathcal{T}/\Z^N}}|K|\int_\R m_K^n(\xi)d\xi
=\frac12\|v_\delta(t_n)\|_{L^2(\T^N)}^2.
\end{equation}
We have used Remark~\ref{rk:suppmnK} to justify the integration by parts in the term with the measure $m^n_K$. The term
\begin{equation}\label{squareFlux}
\sum_{K\in{\mathcal{T}/\Z^N}}\sum_{L\in\mathcal{N}(K)} a_{K\to L}^n(\xi)
\end{equation}
related to the flux term in \eqref{KiFVhalf} has vanished. Indeed, \eqref{squareFlux} is equal to 
\begin{equation}\label{squareFlux2}
\frac12\sum_{K\in{\mathcal{T}/\Z^N}}\sum_{L\in\mathcal{N}(K)} a_{K\to L}^n(\xi)+a_{L\to K}^n(\xi)
\end{equation}
by relabelling of the indexes of summation. All the arguments in \eqref{squareFlux2} cancel individually in virtue of the conservative symmetry property \eqref{conservativesym} of $A_{K\to L}(v,w)$. Indeed, one can check that $a_{K\to L}$ inherits this property, \textit{i.e.}
\begin{equation}\label{conservativesyma}
a_{K\to L}(\xi,v,w)=-a_{L\to K}(\xi,w,v),\quad K,L\in\mathcal{T},\; v,w\in\R,
\end{equation} 
thanks to the explicit formula \eqref{newamon}. To obtain the equation for the balance of energy corresponding to the stochastic forcing, we use the equation
\begin{equation}\label{EulerSDE}
v^{n+1}_K=v^{n+1/2}_K+(\Delta t_n)^{1/2} g_{k,K}(v^n_K)X^{n+1}_k,
\end{equation}
which follows from the equation of the scheme~\eqref{FVscheme} and the definition of $v^{n+1/2}_K$ by \eqref{vhalf}. Taking the square of both sides of \eqref{EulerSDE} and using the independence of $X^{n+1}_k$ and $v^{n+1/2}_i$, we obtain the identity
\begin{equation}\label{eq:FVEnergyEestimate1}
\frac12\E\|v_\delta(t_{n+1})\|_{L^2(\T^N)}^2
=\frac12\E\|v_\delta^\flat(t_{n+1})\|_{L^2(\T^N)}^2
+\frac{\Delta t_n}{2}\E\sum_{K\in{\mathcal{T}/\Z^N}}|K|\sum_{k\geq 1} |g_{k,K}(v^n_K)|^2.
\end{equation}
Adding \eqref{eq:FVEnergyEestimate0} to \eqref{eq:FVEnergyEestimate1} gives \eqref{eq:FVEnergyEestimate}. \qed\bigskip

\begin{remark} Note that \eqref{EulerSDE} also gives
\begin{equation}\label{eq:FVEnergyEestimate11}
\frac12\E\|v_\delta(t_{n+1})-v_\delta^\flat(t_{n+1})\|_{L^2(\T^N)}^2
=
\frac{\Delta t_n}{2}\E\sum_{K\in{\mathcal{T}/\Z^N}}|K|\sum_{k\geq 1} |g_{k,K}(v^n_K)|^2,
\end{equation}
for all $0\leq n\leq N_T$.
\end{remark}
\textbf{Proof of Proposition~\ref{prop:ControlDissipation}.}  
We begin with the proof of the estimates \eqref{eq:ControlSpace} and \eqref{eq:ControlTime}. 
Multiplying Equation~\eqref{KiFVhalf} by $\bar{\mathtt{f}}^n_K:=1-\mathtt{f}^n_K$, we obtain 
\begin{equation}\label{KiFVhalfbar}
|K|\bar{\mathtt{f}}^n_K(\xi)\mathtt{f}^{n+1/2}_K(\xi)+\Delta t_n\bar{\mathtt{f}}^n_K(\xi)\sum_{L\in\mathcal{N}(K)}a_{K\to L}^n(\xi)=|K|\Delta t_n\,\bar{\mathtt{f}}^n_K(\xi)\partial_\xi m^n_K(\xi).
\end{equation}
Next, we multiply \eqref{KiFVhalfbar} by $(\xi-v^n_K)$ and sum the result over $\xi$, $K$. We use the first identity
\begin{equation}\label{fE00}
\int_\R (\xi-v^n_K)\bar{\mathtt{f}}^n_K(\xi)\partial_\xi m^n_K(\xi) d\xi
=\int_\R (\xi-v^n_K)_+\partial_\xi m^n_K(\xi) d\xi
=-\int_\R \bar{\mathtt{f}}^n_K(\xi) m^n_K(\xi) d\xi,
\end{equation}
(once again, we use the fact that $m^n_K$ is compactly supported to do the integration by parts in \eqref{fE00}, \textit{cf.} Remark~\ref{rk:suppmnK}) and the second identity
\begin{equation*}
\int_\R (\xi-v^n_K)\bar{\mathtt{f}}^n_K(\xi)\mathtt{f}^{n+1/2}_K(\xi) d\xi=\frac12 (v^{n+1/2}_K-v^n_K)_+^2,
\end{equation*}
to obtain
\begin{multline}\label{fE1}
\frac12\left\|\big[v_\delta^\flat(t_{n+1})-v_\delta(t_n)\big]_+\right\|_{L^2(\T^N)}^2+\Delta t_n\sum_{K\in\mathcal{T}/\Z^N}|K|\int_\R \bar{\mathtt{f}}^n_K(\xi) m^n_K(\xi) d\xi\\
=-\Delta t_n\sum_{K\in\mathcal{T}/\Z^N}\sum_{L\in\mathcal{N}(K)} \int_\R (\xi-v^n_K)_+ a_{K\to L}^n(\xi) d\xi.
\end{multline}
We transform the right-hand side of \eqref{fE1} by integration by parts in $\xi$: this gives, thanks to \eqref{defPhinum}-\eqref{defPhinumn}, the term
\begin{equation}\label{fE12}
-\Delta t_n\sum_{K\in\mathcal{T}/\Z^N}\sum_{L\in\mathcal{N}(K)} \int_\R \bar{\mathtt{f}}^n_K(\xi) \Phi_{K\to L}^n(\xi) d\xi.
\end{equation}
Then we can relabel the indices in \eqref{fE12} and use the conservative symmetry relation (consequence of \eqref{conservativesyma})
\begin{equation}\label{conservativesymPhi}
\Phi_{K\to L}(\xi,v,w)=-\Phi_{L\to K}(\xi,w,v),
\end{equation}
to see that 
\begin{multline}\label{fE13}
\frac12\left\|\big[v_\delta^\flat(t_{n+1})-v_\delta(t_n)\big]_+\right\|_{L^2(\T^N)}^2+\Delta t_n\sum_{K\in\mathcal{T}/\Z^N}|K|\int_\R \bar{\mathtt{f}}^n_K(\xi) m^n_K(\xi) d\xi\\
=\frac12\Delta t_n\sum_{K\in\mathcal{T}/\Z^N}\sum_{L\in\mathcal{N}(K)} \int_\R 
(\bar{\mathtt{f}}^n_L(\xi)-\bar{\mathtt{f}}^n_K(\xi)) 
\Phi_{K\to L}^n(\xi) d\xi.
\end{multline}
Note that the integrand $(\bar{\mathtt{f}}^n_L(\xi)-\bar{\mathtt{f}}^n_K(\xi)) 
\Phi_{K\to L}^n(\xi) $ is non-negative thanks to the monotony properties of $A_{K\to L}$ and \eqref{eq:numEntropyFlux}.  At this stage, in order to deduce \eqref{eq:ControlSpace} from \eqref{fE13}, we have to prove that, under the CFL condition \eqref{CFLstrong}, a fraction of the right-hand side of \eqref{fE13} controls the term
$$
\frac12\left\|\big[v_\delta^\flat(t_{n+1})-v_\delta(t_n)\big]_+\right\|_{L^2(\T)}^2,
$$
(see the estimate \eqref{EstimTimeWBV2} below). To this end, we integrate Equation~\eqref{KiFVhalfbar} over $\xi\in\R$. This gives
\begin{equation}\label{fE2}
|K|\big[v_K^{n+1/2}-v_K^n\big]_+ + \Delta t_n\sum_{L\in\mathcal{N}(K)}\int_\R \bar{\mathtt{f}}^n_K(\xi) a_{K\to L}^n(\xi)d\xi\leq 0,
\end{equation}
which reads also
\begin{equation*}
|K|\big[v_K^{n+1/2}-v_K^n\big]_+ \leq - \Delta t_n\sum_{L\in\mathcal{N}(K)}\Phi_{K\to L}^n(v^n_K)
\end{equation*}
by \eqref{defPhinumn} (note that it is also the discrete entropy inequality \eqref{eq:discreteEntropyInequality} with $\xi=v^n_K$). Taking the square, using the Cauchy-Schwarz Inequality and summing over $K\in\mathcal{T}/\Z^N$, we deduce that
\begin{equation}\label{EstimTimeWBV}
\left\|\big[v_\delta^\flat(t_{n+1})-v_\delta(t_n)\big]_+\right\|_{L^2(\T^N)}^2
\leq 
\Delta t_n\sum_{K\in\mathcal{T}/\Z^N}\Delta t_n\frac{|\partial K|}{|K|}\sum_{L\in\mathcal{N}(K)}\frac{|\Phi_{K\to L}^n(v^n_K)|^2}{|K|L|}.
\end{equation}
Next, we note that $|\Phi_{K\to L}^n(v^n_K)|^2$ is non-trivial only if $v^n_K<v^n_L$. In that case, it can be decomposed as
\begin{equation}\label{PhiSquare1}
|\Phi_{K\to L}^n(v^n_K)|^2=-2\int_{v^n_K}^{v^n_L} \Phi_{K\to L}^n(\xi)\partial_\xi\Phi_{K\to L}^n(\xi)d\xi
=2\int_{v^n_K}^{v^n_L} \Phi_{K\to L}^n(\xi)a_{K\to L}^n(\xi)d\xi,
\end{equation}
which is bounded by
\begin{equation}\label{PhiSquare2}
2\sup_{\xi\in\R}|a^n_{K\to L}(\xi)|\int_{v^n_K}^{v^n_L} |\Phi_{K\to L}^n(\xi)|d\xi
=2\sup_{\xi\in\R}|a^n_{K\to L}(\xi)|\int_\R (\bar{\mathtt{f}}^n_L-\bar{\mathtt{f}}^n_K) \Phi_{K\to L}^n(\xi) d\xi.
\end{equation}
Under the CFL condition~\eqref{CFLstrong}, the estimate
\eqref{EstimTimeWBV} can be completed into
\begin{multline}\label{EstimTimeWBV2}
\frac12\left\|\big[v_\delta^\flat(t_{n+1})-v_\delta(t_n)\big]_+\right\|_{L^2(\T^N)}^2\\
\leq 
 (1-\theta)\frac12\Delta t_n\sum_{K\in\mathcal{T}/\Z^N}\sum_{L\in\mathcal{N}(K)} \int_\R 
(\bar{\mathtt{f}}^n_L(\xi)-\bar{\mathtt{f}}^n_K(\xi)) 
\Phi_{K\to L}^n(\xi) d\xi.
\end{multline}
Using \eqref{fE13} then, we deduce the two estimates \eqref{eq:ControlSpace}-\eqref{eq:ControlTime}.\medskip

To prove the estimates \eqref{eq:ControlSpacebar} and \eqref{eq:ControlTimebar}, we proceed similarly: we start from the following equation on $\bar{\mathtt{f}}^n_K$, which is equivalent to~\eqref{KiFVhalf}: 
\begin{equation}\label{KiFVhalfbarreal}
|K|(\bar{\mathtt{f}}^{n+1/2}_K(\xi)-\bar{\mathtt{f}}^n_K(\xi))+\Delta t_n\sum_{L\in\mathcal{N}(K)}\bar a_{K\to L}^n(\xi)=-|K|\Delta t_n\, \partial_\xi m^n_K(\xi).
\end{equation}
Then we multiply Eq.~\eqref{KiFVhalfbarreal} by $\mathtt{f}^n_K$, to obtain
\begin{equation}\label{KiFVhalfbarbar}
|K| \mathtt{f}^n_K\bar{\mathtt{f}}^{n+1/2}_K(\xi)+\Delta t_n \sum_{L\in\mathcal{N}(K)}\mathtt{f}^n_K\bar a_{K\to L}^n(\xi)=-|K|\Delta t_n\,\mathtt{f}^n_K \partial_\xi m^n_K(\xi),
\end{equation}
which is the analogue to \eqref{KiFVhalfbar}. In a first step, we multiply \eqref{KiFVhalfbarbar} by $(v^n_K-\xi)$ and sum the result over $\xi\in\R$, $K\in\mathcal{T}/\Z^N$. This gives (compare to \eqref{fE1}-\eqref{fE13})
\begin{align}
\frac12\left\|\big[v_\delta^\flat(t_{n+1})-v_\delta(t_n)\big]_-\right\|_{L^2(\T^N)}^2&+\Delta t_n\sum_{K\in\mathcal{T}/\Z^N}|K|\int_\R \mathtt{f}^n_K(\xi) m^n_K(\xi) d\xi\nonumber\\
=&-\Delta t_n\sum_{K\in\mathcal{T}/\Z^N}\sum_{L\in\mathcal{N}(K)}\mathtt{f}^n_K\overline{\Phi}^n_{K\to L}(\xi)\nonumber\\
=&\frac12\Delta t_n\sum_{K\in\mathcal{T}/\Z^N}\sum_{L\in\mathcal{N}(K)}(\mathtt{f}^n_L-\mathtt{f}^n_K)\overline{\Phi}^n_{K\to L}(\xi).\label{fE1bar}
\end{align}
To conclude to \eqref{eq:ControlSpacebar}-\eqref{eq:ControlTimebar} under the CFL condition \eqref{CFLstrongbar}, we proceed as in \eqref{fE2}-\eqref{EstimTimeWBV2} above, with the minor difference that, instead of the identity $\partial_\xi\overline{\Phi}^n_{K\to L}(\xi)=\bar{a}^n_{K\to L}(\xi)$, we use the formula
$\partial_\xi\overline{\Phi}^n_{K\to L}(\xi)=\bar{b}^n_{K\to L}(\xi,\xi)$ (see \eqref{defbarbnKL}) when we develop $|\overline{\Phi}^n_{K\to L}(v^n_K)|^2$.  \qed\bigskip


\begin{remark} A slight modification of the lines \eqref{PhiSquare1}-\eqref{PhiSquare2} in the proof above shows that
\begin{equation}\label{PhiSquare3}
|\Phi_{K\to L}^n(\xi\vee v^n_K)|^2
\leq 2\sup_{\xi\in\R}|a^n_{K\to L}(\xi)|\int_\R (\bar{\mathtt{f}}^n_L-\bar{\mathtt{f}}^n_K) \Phi_{K\to L}^n(\xi) d\xi,
\end{equation}
for all $\xi\in\R$. This estimate will be used in the proof of Lemma~\ref{lem:consisSPACE} below.
\label{rk:PhiSquare}\end{remark}

\textbf{Proof of Corollary~\ref{cor:EnergyEstimate}.}  Assume that \eqref{CFLstrongU} is satisfied. It is clear, in virtue of the estimate \eqref{alphaKpK} and the bound \eqref{BOUNDa} and \eqref{BOUNDastar} on $a^n_{K\to L}$ and $\bar{b}^n_{K\to L}$, that \eqref{CFLstrongU} implies the CFL conditions \eqref{CFLstrong} and \eqref{CFLstrongbar}. Besides, due to \eqref{D0num}, we have the bound
\begin{equation}\label{D0numnum}
\sum_{K\in\mathcal{T}/\Z^N}|K|\sum_{k\geq 1} |g_{k,K}(v^n_K)|^2\leq D_0.
\end{equation}
This gives 
$$
\sum_{n=0}^{N_T-1}\Delta t_n\sum_{K\in\mathcal{T}/\Z^N}|K|\sum_{k\geq 1} |g_{k,K}(v^n_K)|^2\leq D_0 T,
$$
which, inserted in the energy estimate \eqref{eq:FVEnergyEestimate}, shows that
$$
\E\mathcal{E}(T)\leq \frac{1}{2}\|v_\delta(0)\|_{L^2(\T^N)}^2+\frac12 D_0 T.
$$
By addition of the estimates \eqref{eq:ControlSpace}-\eqref{eq:ControlSpacebar} and \eqref{eq:ControlTime}-\eqref{eq:ControlTimebar} respectively, we obtain therefore \eqref{eq:ControlSpaceAdd} and \eqref{eq:ControlTimeAdd}. There remains to prove \eqref{eq:ControlTimeAdd01}. For that purpose, we use \eqref{eq:FVEnergyEestimate11} and \eqref{D0numnum} to obtain  
\begin{equation}\label{endcor}
\E\|v_\delta(t_{n+1})-v_\delta^\flat(t_{n+1})\|_{L^2(\T^N)}^2=\E\|v_\delta(t_{n+1})\|_{L^2(\T^N)}^2-\E\|v_\delta^\flat(t_{n+1})\|_{L^2(\T^N)}^2\leq D_0\Delta t_n.
\end{equation}
Summing \eqref{endcor} over $0\leq n<N_T$ and using \eqref{eq:ControlTimeAdd} yields \eqref{eq:ControlTimeAdd01}. \qed

\section{Approximate kinetic equation}\label{sec:appki}

\subsection{Discrete unknown}\label{sec:discretef}

In this section we will show that, when $\delta\to0$, some discrete kinetic unknowns $f_\delta$ associated to the scheme \eqref{FVscheme} are approximate kinetic solutions. There may be several way to define $f_\delta$: it depends for example on the manner in which the discrete data $f^n_K$ are assembled by interpolation. One of the constraints due to our definition of approximate kinetic solution is the formulation ``at fixed $t$" \eqref{eq:kineticfpreappt} (in opposition to a formulation which would be weak in time). To establish such a formulation in our context, a minimal amount of regularity of the function $t\mapsto f_\delta(t)$ is required (in particular, for all $\varphi$, $\<f_\delta(t),\varphi\>$ should be a c{\`a}dl{\`a}g process). For $t\in[t_n,t_{n+1})$, we will therefore consider the function $f_\delta(t)$ defined as the interpolation
\begin{equation}\label{deffh}
f_\delta(x,t,\xi)=\frac{t-t_n}{\Delta t_n}\mathbf{1}_{v^\sharp_\delta(x,t)>\xi}+\frac{t_{n+1}-t}{\Delta t_n}\mathbf{1}_{v_\delta(x,t)>\xi},\quad \xi\in\R, x\in\T^N,
\end{equation}
where $v^\sharp_\delta(x,t)$ is given by
\begin{equation}\label{def:barvni}
v^\sharp_\delta(x,t)=v^{n+1/2}_K+g_{k,K}(v^n_K)(\beta_k(t)-\beta_k(t_n)),\quad t_n\leq t< t_{n+1}, x\in K.
\end{equation}
We set $v^\sharp_K(t)=v^\sharp_\delta(x,t)$, $x\in K$. Then, for $t\in[t_n,t_{n+1})$ , $t\mapsto v^\sharp_K(t)$ is itself an interpolation between $v^{n+1/2}_K$ and $v^{n+1}_K$. We also denote by $\mathtt{f}_\delta$ the piecewise constant function
\begin{equation}\label{moredelta}
\mathtt{f}_\delta(x,t,\xi)=\mathtt{f}^n_K=\mathbf{1}_{v_\delta(x,t)>\xi},\quad x\in K, t\in[t_n,t_{n+1}).
\end{equation}
We check first that $f_\delta$ and $\mathtt{f}_\delta$ are close to each other.

\begin{lemma} Let $u_0\in L^\infty(\T^N)$, $T>0$. Assume that \eqref{D1}, \eqref{multiplicativeNoise}, \eqref{D0plus}, \eqref{AALip}, and \eqref{CFLstrongU} are satisfied. For $\delta\in \mathfrak{d}_T$, assume \eqref{alphaK} and \eqref{alphapK}. Let $(v_\delta(t))$ be the numerical unknown defined by \eqref{FVscheme}-\eqref{FVIC}-\eqref{defvh} and let $f_\delta$, $\mathtt{f}_\delta$ be defined by \eqref{deffh}-\eqref{moredelta}. Then 
\begin{multline}\label{eq:fdeltafdelta}
\E\int_0^T\int_{\T^N}\left|\int_\R|f_\delta(x,t,\xi)-\mathtt{f}_\delta(x,t,\xi)|d\xi\right|^2 dx dt\\
\leq
\left[\theta^{-1}\|v_\delta(0)\|_{L^2(\T^N)}^2+D_0 T(1+\theta^{-1})\right]\left[\sup_{0\leq n<N_T}\Delta t_n\right].
\end{multline}
\label{lem:fdeltafdelta}
\end{lemma}

\textbf{Proof of Lemma~\ref{lem:fdeltafdelta}.} Since 
$$
f_\delta(t)-\mathtt{f}_\delta(t)=\frac{t-t_n}{\Delta t_n}(\mathbf{1}_{v^\sharp_\delta(t)>\xi}-\mathbf{1}_{v_\delta(t)>\xi}),
$$
for $t\in[t_n,t_{n+1})$ and since the factor $\frac{t-t_n}{\Delta t_n}$ is less than $1$, the quantity we want to estimate is bounded by the following $L^2$-norm:
\begin{equation}\label{fdeltaL2L2}
\E\int_0^T\int_{\T^N}\left|\int_\R|f_\delta(x,t,\xi)-\mathtt{f}_\delta(x,t,\xi)|d\xi\right|^2 dx dt
\leq\E\int_0^T\|v^\sharp_\delta(t)-v_\delta(t)\|_{L^2(\T^N)}^2 dt.
\end{equation}
By definition of $v^\sharp_\delta(t)$ and independence and \eqref{D0num}, we obtain
\begin{multline*}
\E\int_0^T\int_{\T^N}\left|\int_\R|f_\delta(x,t,\xi)-\mathtt{f}_\delta(x,t,\xi)|d\xi\right|^2 dx dt
\leq D_0\sum_{n=0}^{N_T-1}\int_{t_n}^{t_{n+1}}|t-t_n| dt\\
+\E\sum_{n=0}^{N_T-1}\Delta t_n\left\|v_\delta^\flat(t_{n+1})-v_\delta(t_n)\right\|_{L^2(\T^N)}^2.
\end{multline*}
Using the temporal estimate \eqref{eq:ControlTimeAdd}, we deduce \eqref{eq:fdeltafdelta}. \qed\bigskip

\begin{remark} Note for a future use (\textit{cf.} \eqref{epsSTO26}) that we have just proved the estimate
\begin{equation}\label{timeBVvbarv}
\E\int_0^T\|v^\sharp_\delta(t)-v_\delta(t)\|_{L^2(\T^N)}^2 dt
\leq 
\left[\theta^{-1}\|v_\delta(0)\|_{L^2(\T^N)}^2+D_0 T(1+\theta^{-1})\right]\left[\sup_{0\leq n<N_T}\Delta t_n\right]
\end{equation}
\end{remark}
To $f_\delta$ we will associate the Young measure
\begin{equation}\label{defnuh}
\nu^\delta_{x,t}(\xi):=-\partial_\xi f_\delta(x,t,\xi)=\frac{t-t_n}{\Delta t_n}\delta(\xi=v^\sharp_\delta(x,t))+\frac{t_{n+1}-t}{\Delta t_n}\delta(\xi=v_\delta(x,t)),
\end{equation}
We also denote by $m_\delta$ the discrete random measure given by
\begin{equation}\label{defmh}
d m_\delta(x,t,\xi)=\sum_{n=0}^{N_T-1}\sum_{K\in\mathcal{T}}\mathbf{1}_{K\times[t_n,t_{n+1})}(x,t)m^n_K(\xi)\,dx dt d\xi.
\end{equation}
Recall the definitions \eqref{ADMdelta1}-\eqref{ADMdelta3} (definition of the set of mesh parameter $\mathfrak{d}_T$ in par\-ti\-cu\-lar), that we will use in all the section.

\begin{proposition}[Discrete kinetic equation] Let $u_0\in L^\infty(\T^N)$, $T>0$. Assume that \eqref{D1}, \eqref{multiplicativeNoise}, \eqref{D0plus}, \eqref{AALip} and \eqref{CFLstrongU} are satisfied. For $\delta\in \mathfrak{d}_T$, assume \eqref{alphaK} and \eqref{alphapK}. Let $(v_\delta(t))$ be the numerical unknown defined by \eqref{FVscheme}-\eqref{FVIC}-\eqref{defvh} and let $f_\delta$, $\nu^\delta$, $m_\delta$ be defined by \eqref{deffh}, \eqref{defnuh}, \eqref{defmh} respectively. Then $f_\delta$ satisfies the following discrete kinetic formulation: for all $t\in[t_n,t_{n+1}]$, $x\in K$, for all $\psi\in C^1_c(\R)$,
\begin{align}
&\<f_\delta(x,t),\psi\>-\<f_\delta(x,t_n),\psi\>\nonumber\\
=&-\frac{1}{|K|}\int_{t_n}^t\int_\R \sum_{L\in\mathcal{N}(K)}a^n_{K\to L}(\xi)\psi(\xi) d\xi ds
-\int_{t_n}^t\int_\R \partial_\xi\psi(\xi) m^n_K(\xi) d\xi ds\nonumber\\
&+\frac{t-t_n}{\Delta t_n}\int_{t_n}^t g_{k,K}(v^n_K)\psi(v^\sharp_\delta(x,s))d\beta_k(s)
+\frac12\frac{t-t_n}{\Delta t_n}\int_{t_n}^t \GG^2_K(v^n_K)\partial_\xi\psi(v^\sharp_\delta(x,s))ds.\label{eq:discretekieq}
\end{align}
\label{prop:discretekieq}\end{proposition}

In \eqref{eq:discretekieq}, $\<f_\delta(x,t),\psi\>$ stands for the product
$$
\int_\R f_\delta(x,t,\xi)\psi(\xi)d\xi.
$$

\textbf{Proof of Proposition~\ref{prop:discretekieq}.} Let $\Psi$ be a primitive for $\psi$ and let $x\in K$, $t\in[t_n,t_{n+1})$. By definition of $f_\delta$, see Equation~\eqref{deffh}, we have
$$
\<f_\delta(x,t),\psi\>-\<f_\delta(x,t_n),\psi\>=\frac{t-t_n}{\Delta t_n}\left[\Psi(v^\sharp(x,t))-\Psi(v_\delta(x,t))\right],
$$
which we decompose as the sum of two terms:
\begin{equation}\label{discretekieq1}
\frac{t-t_n}{\Delta t_n}\left[\Psi(v^\sharp(x,t))-\Psi(v^\flat(x,t_{n+1}))\right],
\end{equation}
and
\begin{equation}\label{discretekieq2}
\frac{t-t_n}{\Delta t_n}\left[\Psi(v^\flat(x,t_{n+1}))-\Psi(v_\delta(x,t))\right].
\end{equation}
We use first the deterministic kinetic formulation \eqref{KiFVhalf}, which we multiply by $\psi(\xi)$. By integration over $\xi\in\R$, it gives
\begin{equation}\label{discretekieq3}
\eqref{discretekieq2}=-\frac{1}{|K|}\int_{t_n}^t\int_\R \sum_{L\in\mathcal{N}(K)}a^n_{K\to L}(\xi)\psi(\xi) d\xi ds
-\int_{t_n}^t\int_\R \partial_\xi\psi(\xi) m^n_K(\xi) d\xi ds.
\end{equation} 
By It\={o}'s Formula on the other hand (\textit{cf.} \eqref{def:barvni}), the term \eqref{discretekieq1} is equal to
\begin{equation}\label{discretekieq4}
\frac{t-t_n}{\Delta t_n}\int_{t_n}^t g_{k,K}(v^n_K)\psi(v^\sharp_\delta(x,s))d\beta_k(s)
+\frac12\frac{t-t_n}{\Delta t_n}\int_{t_n}^t \GG^2_K(v^n_K)\partial_\xi\psi(v^\sharp_\delta(x,s))ds.
\end{equation} 
Summing \eqref{discretekieq3} and \eqref{discretekieq4}, we obtain \eqref{eq:discretekieq}. \qed\bigskip


We will prove now that the Finite Volume scheme \eqref{FVscheme} is consistent with \eqref{stoSCL}. Indeed, we will show, thanks to the estimates obtained in Section~\ref{secEnergyEstimates}, that an approximate kinetic equation for $f_\delta$ in the sense of \eqref{eq:kineticfpreappt} can be deduced from the discrete kinetic formulation \eqref{eq:discretekieq}.

\begin{proposition}[Approximate kinetic equation] Let $u_0\in L^\infty(\T^N)$, $T>0$. Assume that \eqref{D1}, \eqref{multiplicativeNoise}, \eqref{D0plus}, \eqref{AALip} and \eqref{CFLstrongU} are satisfied. For $\delta\in \mathfrak{d}_T$, assume \eqref{alphaK} and \eqref{alphapK}. Let $(v_\delta(t))$ be the numerical unknown defined by \eqref{FVscheme}-\eqref{FVIC}-\eqref{defvh} and let $f_\delta$, $\nu^\delta$, $m_\delta$ be defined by \eqref{deffh}, \eqref{defnuh}, \eqref{defmh} respectively. If $(\delta_m)$ is a sequence in $\mathfrak{d}_T$ which tends to zero according to \eqref{ADMdelta3}, then $(f_{\delta_m})$ is a sequence of approximate solutions to \eqref{stoSCL} of order $\order=2$. Besides, $(f_{\delta_m}(0))$ converges to the equilibrium function $\mathtt{f}_0=\mathbf{1}_{u_0>\xi}$ in $L^\infty(\T^N\times\R)$-weak-*. 
\label{prop:consis}\end{proposition}

\textbf{Proof of Proposition~\ref{prop:consis}.} The last assertion is clear: $(f_{\delta_m}(0))$ converges to the equilibrium function $\mathtt{f}_0=\mathbf{1}_{u_0>\xi}$ in $L^\infty(\T^N\times\R)$-weak-* since, by \eqref{FVIC}, $v^{\delta_m}(0)\to u_0$ a.e. on $\T^N$. We will show that, for all $t\in[0,T]$, for all $\varphi\in C^\infty_c(\T^N\times\R)$,
\begin{align}
\<f_\delta(t),\varphi\>=&\<f_\delta(0),\varphi\>
-\iiint_{\T^N\times[0,t]\times\R}\partial_\xi\varphi(x,\xi)d m_\delta(x,s,\xi)+\eps^\delta(t,\varphi)\nonumber\\
&+\int_0^t \<f_\delta(s),a(\xi)\nabla_x\varphi\>ds\label{ExactSpace}\\
&+\int_0^t\int_{\T^N}\int_\R g_k(x,\xi)\varphi(x,\xi)d\nu^\delta_{x,s}(\xi) dxd\beta_k(s) \label{ExactNoise}\\
&+\frac{1}{2}\int_0^t\int_{\T^N}\int_\R \GG^2(x,\xi)\partial_\xi\varphi(x,\xi)d\nu^\delta_{x,s}(\xi) dx ds,
\label{ExactCorrection}
\end{align}
where the error term $\eps^\delta(t,\varphi)$ satisfies
\begin{equation}\label{pdeltato0}
\lim\limits_{\delta\to 0}\E\left[\sup_{t\in[0,T]}|\eps^\delta(t,\varphi)|^2\right]=0,
\end{equation} 
for all $\varphi\in C^2_c(\T^N\times\R)$. Note that the convergence in probability \eqref{epsto0} follows from \eqref{pdeltato0}. Given $\varphi\in C^2_c(\T^N\times\R)$, we introduce the averages over the cells $K\in\mathcal{T}$
\begin{equation}\label{defphini}
\varphi_K(\xi)=\frac{1}{|K|}\int_{|K|}\varphi(x,\xi)dx,\quad\xi\in\R.
\end{equation}
To prove \eqref{ExactCorrection}, we apply the discrete kinetic equation \eqref{eq:discretekieq} to $\xi\mapsto\varphi(x,\xi)$ for a fixed $x\in K$. Then we sum the result over $x\in\T^N$. By the telescopic formula
$$
\<f_\delta(x,t),\varphi\>-\<f_\delta(x,0),\varphi\>=\sum_{n=0}^{N_T-1} \<f_\delta(x,t\wedge t_{n+1}),\varphi\>-\<f_\delta(x,t\wedge t_n),\varphi\>,
$$
we obtain
\begin{align}
\<f_\delta(t),\varphi\>=&\<f_\delta(x,0),\varphi\>-\iiint_{\T^N\times[0,t]\times\R} \partial_\xi\varphi(\xi) dm_\delta(x,s,\xi) \nonumber\\
&-\sum_{n=0}^{N_T-1}\sum_{K\in\mathcal{T}/\Z^N}\int_{t\wedge t_n}^{t\wedge t_{n+1}}\int_\R\sum_{L\in\mathcal{N}(K)} a^n_{K\to L}(\xi)\varphi_K(\xi) ds d\xi\label{AppSpace}\\
&+\int_0^t\int_{\T^N}\int_{\R\times\R} g_{k,\delta}(x,\xi)\varphi(x,\zeta)d\mu^\delta_{x,s,t}(\xi,\zeta)d\beta_k(s)\label{AppNoise}\\
&+\frac12\int_0^t\int_{\T^N}\int_{\R\times\R} \GG_\delta^2(x,\xi)\partial_\xi\varphi(x,\zeta) d\mu^\delta_{x,s,t}(\xi,\zeta) ds,\label{AppCorrection}
\end{align}
where the measure $\mu^\delta_{x,s,t}$ on $\R\times\R$ is defined by
\begin{equation}\label{defmudouble}
\<\mu^\delta_{x,s,t},\psi\>=\sum_{n=0}^{N_T-1}\frac{t\wedge t_{n+1}-t\wedge t_{n}}{\Delta t_n}
\mathbf{1}_{[t_n,t_{n+1})}(s)
\psi(v_\delta(x,s),v^\sharp_\delta(x,s)),\quad \psi\in C_b(\R^2),
\end{equation}
and the discrete coefficient $g_{k,\delta}(x,\xi)$ is equal to $g_{k,K}(\xi)$ (\textit{cf.} \eqref{defgnum}) when $x\in K$ (similarly, $\GG_\delta(x,\xi):=\GG_K(\xi)$, $x\in K$). Note that $\mu^\delta_{x,s,t}$ is simply the Dirac mass at $(v_\delta(x,s),v^\sharp_\delta(x,s))$, except when $t_l\leq s\leq t$ (where $l$ is the index such that $t_l\leq t< t_{l+1}$), in which case it is the same Dirac mass, with an additional multiplicative factor $\frac{t-t_l}{\Delta t_l}$. \medskip

The term \eqref{AppSpace} is a discrete space derivative: we will show that it is an approximation of the term \eqref{ExactSpace}. The two terms \eqref{AppNoise} and \eqref{AppCorrection} are close to \eqref{ExactNoise} and \eqref{ExactCorrection} respectively. We analyse those terms separately (see Section~\ref{sec:appspace}, Section~\ref{sec:appSto}). The conclusion of the proof of Proposition~\ref{prop:consis} is given in Section~\ref{sec:appcl}.

\subsection{Space consistency}\label{sec:appspace}

\begin{lemma} Let $u_0\in L^\infty(\T^N)$, $T>0$ and $\delta\in \mathfrak{d}_T$. Assume that \eqref{multiplicativeNoise}, \eqref{D0plus}, \eqref{AALip}, \eqref{alphaK}, \eqref{alphapK} and \eqref{CFLstrongU} are satisfied. Then, for all $\varphi\in C^2_c(\T^N\times\R)$, we have
\begin{align}
\int_{\T^N}&\int_0^t\int_\R a(\xi)\cdot\nabla_x\varphi f_\delta(s) dx ds d\xi\nonumber\\
=&-\sum_{n=0}^{N_T-1}\sum_{K\in\mathcal{T}/\Z^N}\int_{t\wedge t_n}^{t\wedge t_{n+1}}\int_\R\sum_{L\in\mathcal{N}(K)} a^n_{K\to L}(\xi)\varphi_K(\xi) d\xi+\eps^\delta_{\mathrm{space},0}(t,\varphi)+\eps^\delta_{\mathrm{space},1}(t,\varphi),\label{consistSPACE}
\end{align}
for all $t\in[0,T]$, with the estimates
\begin{align}
\E\sup_{t\in[0,T]}&|\eps^\delta_{\mathrm{space},0}(t,\varphi)|^2\nonumber\\
&\leq T|{L_A}|^2\|\nabla_x\varphi\|_{L^\infty_{x,\xi}}^2\left[\frac{1}{\theta}\|v_\delta(0)\|_{L^2(\T^N)}^2+D_0 T\left(1+\frac{1}{\theta}\right)\right]\sup_{0\leq n<N_T}\Delta t_n,
\label{eps00SPACE}
\end{align}
and, for all compact $\Lambda\subset\R$, for all $\varphi\in C^\infty_c(\T^N\times\R)$ supported in $\T^N\times \Lambda$,
\begin{equation}
\E\sup_{t\in[0,T]}|\eps^\delta_{\mathrm{space},1}(t,\varphi)|^2
\leq \frac{16{L_A}T}{\alpha_N^2}|\Lambda|^2\|\partial_{\xi}\nabla_x\varphi\|_{L^\infty_{x,\xi}}^2\left[\frac{1}{\theta}\|v_\delta(0)\|_{L^2(\T^N)}^2+\frac{2 D_0 T}{\theta}\right]\, h.
\label{eps11SPACE}
\end{equation}
\label{lem:consisSPACE}\end{lemma}

\textbf{Proof of Lemma~\ref{lem:consisSPACE}.} To begin with, we replace $f_\delta$ by $\mathtt{f}_\delta$ in the left-hand side of \eqref{consistSPACE}. This accounts for the first error term 
$$
\eps^\delta_{\mathrm{space},0}(t,\varphi)=\int_{\T^N}\int_0^t\int_\R a(\xi)\cdot\nabla_x\varphi (f_\delta(s)-\mathtt{f}_\delta(s)) dx ds d\xi.
$$
Thanks to Lemma~\ref{lem:fdeltafdelta}, we have the estimate \eqref{eps00SPACE} for $\eps^\delta_{\mathrm{space},0}(t,\varphi)$. Then, we use the following development:
\begin{align}
\int_{\T^N}\int_0^t\int_\R a(\xi)&\cdot\nabla_x\varphi \mathtt{f}_\delta(s)) dx ds d\xi\nonumber\\
=&\sum_{n=0}^{N_T-1}\int_{t\wedge t_n}^{t\wedge t_{n+1}}\int_\R \sum_{K\in\mathcal{N}(K)}\int_K a(\xi)\cdot\nabla_x\varphi \mathtt{f}_\delta(s) dx ds d\xi,\label{devSpace00}
\end{align}
Since $\mathtt{f}_\delta(s)$ has a constant value $\mathtt{f}^n_K$ in $K\times[t_n,t_{n+1})$, we obtain, thanks to the Stokes formula,
\begin{align}
\int_{\T^N}\int_0^t\int_\R a(\xi)&\cdot\nabla_x\varphi\mathtt{f}_\delta(s) dx ds d\xi\nonumber\\
=&\sum_{n=0}^{N_T-1} \int_{t\wedge t_n}^{t\wedge t_{n+1}}\int_\R\sum_{K\in\mathcal{N}(K)}\sum_{L\in\mathcal{N}(K)} a^*_{K\to L}(\xi)\varphi_{K|L} \mathtt{f}^n_K ds d\xi,\label{devSpace01}
\end{align}
where $a^*_{K\to L}(\xi)$ is defined by \eqref{astarKL} and $\varphi_{K|L}$ is the mean-value of $\varphi$ over $K|L$:
$$
\varphi_{K|L}(\xi)=\frac{1}{|K|L|}\int_{K|L}\varphi(x,\xi)d\mathcal{H}^{N-1}(x).
$$
We add a corrective term to \eqref{devSpace01} to obtain
\begin{align}
\int_{\T^N}\int_0^t\int_\R a(\xi)&\cdot\nabla_x\varphi\mathtt{f}_\delta(s) dx ds d\xi\nonumber\\
=&\sum_{n=0}^{N_T-1} \int_{t\wedge t_n}^{t\wedge t_{n+1}}\int_\R\sum_{K\in\mathcal{N}(K)}\sum_{L\in\mathcal{N}(K)} a^*_{K\to L}(\xi)(\varphi_{K|L}-\varphi_{K}) \mathtt{f}^n_K ds d\xi.\label{devSpace02}
\end{align}
Equation~\eqref{devSpace02} follows indeed from \eqref{devSpace01} by the anti-symmetry property \eqref{conservativesyma} of $a_{K\to L}$. Note that Equation~\eqref{devSpace02} is more natural than Equation~\eqref{devSpace01} (when one has in mind the decomposition of a volume integral over each cells $K$), by use of the correspondence 
$$
a(\xi)\cdot\nabla_x\varphi\simeq \sum_{L\in\mathcal{N}(K)} a^*_{K\to L}(\xi)(\varphi_{K|L}-\varphi_{K})\mbox{ in }K.
$$
By \eqref{conservativesyma}, the discrete convective term in \eqref{AppSpace} is
\begin{equation}\label{AppSpace2}
\sum_{n=0}^{N_T-1}\sum_{K\in\mathcal{T}/\Z^N}\int_{t\wedge t_n}^{t\wedge t_{n+1}}\int_\R\sum_{L\in\mathcal{N}(K)} a^n_{K\to L}(\xi)(\varphi_{K|L}-\varphi_K) ds d\xi.
\end{equation}
To estimate how close is the right-hand side of \eqref{devSpace02} to \eqref{AppSpace2}, we have to compare $a^n_{K\to L}(\xi)$ and $a^*_{K\to L}(\xi)\mathtt{f}^n_K(\xi)$. Let $\gamma\in W^{1,1}(\R_\xi)$. If $v^n_K\leq v^n_L$, then
\begin{equation}
\int_\R\gamma(\xi)\left[a^*_{K\to L}(\xi)\mathtt{f}^n_K(\xi)-a^n_{K\to L}(\xi)\right]d\xi
=-\int_\R\gamma(\xi)\bar{\mathtt{f}}^n_K(\xi) a^n_{K\to L}(\xi)d\xi
\end{equation}
by the consistency hypothesis \eqref{consistencya2} and the support condition~\eqref{supportaKLxi}. Using an integration by parts and \eqref{defPhinum}, we obtain
\begin{equation}\label{gammaKL}
\int_\R\gamma(\xi)\left[a^*_{K\to L}(\xi)\mathtt{f}^n_K(\xi)-a^n_{K\to L}(\xi)\right]d\xi
=\int_\R\gamma'(\xi)\Phi^n_{K\to L}(\xi\vee v^n_K)d\xi.
\end{equation}
Similarly, if $v^n_L\leq v^n_K$, then
\begin{equation}\label{gammaLK}
\int_\R\gamma(\xi)\left[a^*_{K\to L}(\xi)\mathtt{f}^n_K(\xi)-a^n_{K\to L}(\xi)\right]d\xi
=-\int_\R\gamma'(\xi)\overline{\Phi}^n_{K\to L}(\xi\wedge v^n_K)d\xi.
\end{equation}
We deduce that \eqref{consistSPACE} is satisfied with an error term
\begin{align}
\eps^\delta_{\mathrm{space},1}(t,\varphi)
:=\int_{\T^N}\int_0^t&\int_\R a(\xi)\cdot\nabla_x\varphi\, \mathtt{f}_\delta(s) dx ds d\xi\nonumber\\
&+\sum_{n=0}^{N_T-1}\sum_{K\in\mathcal{T}/\Z^N}\int_{t\wedge t_n}^{t\wedge t_{n+1}}\int_\R\sum_{L\in\mathcal{N}(K)} a^n_{K\to L}(\xi)\varphi_K(\xi) d\xi,\label{defepsspace}
\end{align}
which is bounded as follows:
\begin{align}
|\eps^\delta_{\mathrm{space},1}(t,\varphi)|
\leq
\sum_{n=0}^{N_T-1}\int_{t\wedge t_n}^{t\wedge t_{n+1}}&\sum_{K\in\mathcal{T}/\Z^N}\int_\R\sum_{L\in\mathcal{N}(K)} |\partial_\xi\varphi_{K|L}-\partial_\xi\varphi_K| \nonumber\\
\times &\Big[\mathbf{1}_{v^n_K\leq v^n_L}|\Phi^n_{K\to L}(\xi\vee v^n_K)|
+\mathbf{1}_{v^n_L< v^n_K}|\overline{\Phi}^n_{K\to L}(\xi\wedge v^n_K)| \Big] d\xi,\label{boundepsspace1}
\end{align}
By \eqref{hsizemesh}, we have
$$
|\partial_\xi\varphi_{K|L}(\xi)-\partial_\xi\varphi_K(\xi)|\leq \|\partial_\xi\nabla_x\varphi(\cdot,\xi)\|_{L^\infty_{x}}h,
$$
for all $\xi\in\R$. If $\varphi$ is compactly supported in $\T^N\times\Lambda$, we obtain thus the bound 
\begin{equation}\label{boundepsspace2}
|\eps^\delta_{\mathrm{space},1}(t,\varphi)|
\leq
\|\partial_\xi\nabla_x\varphi\|_{L^\infty_{x,\xi}}|\Lambda| B_\mathrm{space}\, h,
\end{equation}
where $B_\mathrm{space}$ is equal to
$$
\sum_{n=0}^{N_T-1}\Delta t_n\sum_{K\in\mathcal{T}/\Z^N}\sum_{L\in\mathcal{N}(K)} 
\Big[\mathbf{1}_{v^n_K\leq v^n_L}\sup_{\xi\in\R}|\Phi^n_{K\to L}(\xi\vee v^n_K)|
+\mathbf{1}_{v^n_L< v^n_K}\sup_{\xi\in\R}|\overline{\Phi}^n_{K\to L}(\xi\wedge v^n_K)| \Big].
$$
We seek for a bound of order $h^{-1/2}$ on $B_\mathrm{space}$. For notational convenience we will estimate only the first part
$$
B_\mathrm{space}^1:=\sum_{n=0}^{N_T-1}\Delta t_n\sum_{K\in\mathcal{T}/\Z^N}\sum_{L\in\mathcal{N}(K)} 
\mathbf{1}_{v^n_K\leq v^n_L}\sup_{\xi\in\R}|\Phi^n_{K\to L}(\xi\vee v^n_K)|,
$$
since the bound on the second part in $B_\mathrm{space}$ will be similar. By the Cauchy Schwarz inequality, we have
\begin{align*}
\E|B_\mathrm{space}^1|^2
\leq \sum_{n=0}^{N_T-1}\Delta t_n&\sum_{K\in\mathcal{T}/\Z^N}|\partial K|\\
&\times \E\sum_{n=0}^{N_T-1}\Delta t_n\sum_{K\in\mathcal{T}/\Z^N}\sum_{L\in\mathcal{N}(K)} 
\frac{\mathbf{1}_{v^n_K\leq v^n_L}}{|K|L|}\sup_{\xi\in\R}|\Phi^n_{K\to L}(\xi\vee v^n_K)|^2.
\end{align*}
We use the estimate \eqref{PhiSquare3}, which gives 
$$
|\Phi_{K\to L}^n(\xi\vee v^n_K)|^2
\leq 2{L_A} |K|L| \int_\R (\bar{\mathtt{f}}^n_L-\bar{\mathtt{f}}^n_K) \Phi_{K\to L}^n(\xi) d\xi,
$$
thanks to \eqref{BOUNDa}. We also use \eqref{alphaKpK}, and get
$$
\E|B_\mathrm{space}^1|^2
\leq \frac{2{L_A}T}{\alpha_N^2 h} \E\sum_{n=0}^{N_T-1}\Delta t_n\sum_{K\in\mathcal{T}/\Z^N}\sum_{L\in\mathcal{N}(K)} \int_\R (\bar{\mathtt{f}}^n_L-\bar{\mathtt{f}}^n_K) \Phi_{K\to L}^n(\xi) d\xi.
$$
With \eqref{eq:ControlSpaceAdd} and \eqref{boundepsspace2}, we conclude to \eqref{eps11SPACE}. \qed

\subsection{Stochastic terms}\label{sec:appSto}
\begin{lemma} Let $u_0\in L^\infty(\T^N)$, $T>0$ and $\delta\in \mathfrak{d}_T$. Assume that \eqref{multiplicativeNoise}, \eqref{D0plus}, \eqref{AALip} and \eqref{CFLstrongU} are satisfied. Then, for all $\varphi\in C^2_c(\T^N\times\R)$, we have
\begin{align}
&\int_0^t\int_{\T^N}\int_{\R\times\R} g_{k,\delta}(x,\xi)\varphi(x,\zeta)d\mu^\delta_{x,s,t}(\xi,\zeta)d\beta_k(s)\nonumber\\
=&\int_0^t\int_{\T^N}\int_\R g_k(x,\xi)\varphi(x,\xi)d\nu^\delta_{x,s}(\xi) dxd\beta_k(s)+\eps^\delta_{\mathrm{W},1}(t,\varphi)+\eps^\delta_{\mathrm{W},2}(t,\varphi),
\label{consistSTONoise}
\end{align}
and
\begin{align}
&\int_0^t\int_{\T^N}\int_{\R\times\R} \GG_\delta^2(x,\xi)\partial_\xi\varphi(x,\zeta) d\mu^\delta_{x,s,t}(\xi,\zeta) ds\nonumber\\
=&\int_0^t\int_{\T^N}\int_\R \GG^2(x,\xi)\partial_\xi\varphi(x,\xi)d\nu^\delta_{x,s}(\xi) dx ds+\eps^\delta_{\mathrm{W},3}(t,\varphi)+\eps^\delta_{\mathrm{W},4}(t,\varphi),
\label{consistSTOCorrection}
\end{align}
where 
\begin{align}
\E\left[\sup_{t\in[0,T]}|\eps^\delta_{\mathrm{W},1}(t,\varphi)|^2\right]\leq&
4 D_1 T \|\varphi\|_{L^\infty_{x,\xi}}^2 h^2
+2D_0 \|\varphi\|_{L^\infty_{x,\xi}}^2\left[\sup_{0\leq n<N_T}\Delta t_n\right],\label{epsW11}\\
\E\left[\sup_{t\in[0,T]}|\eps^\delta_{\mathrm{W},3}(t,\varphi)|^2\right]\leq&
4 D_1 T \|\partial_\xi\varphi\|_{L^\infty_{x,\xi}}^2 h^2
+2D_0 \|\partial_\xi\varphi\|_{L^\infty_{x,\xi}}^2\left[\sup_{0\leq n<N_T}\Delta t_n\right],\label{epsW33}
\end{align}
and
\begin{align}
\E\left[\sup_{t\in[0,T]}|\eps^\delta_{\mathrm{W},2}(t,\varphi)|^2\right]\leq&
2D_0 \|\varphi\|_{L^\infty_{x,\xi}}^2\left[\sup_{0\leq n<N_T}\Delta t_n\right]\nonumber\\
&+8\left\{D_1\|\varphi\|_{L^\infty_{x,v}}^2+D_0\|\partial_\xi\varphi\|_{L^\infty_{x,\xi}}^2\right\}\nonumber\\
&\times\left[\frac{1}{\theta}\|v_\delta(0)\|_{L^2(\T)}^2+\frac{3 D_0 T}{\theta}\right]\left[\sup_{0\leq n<N_T}\Delta t_n\right]^{1/2}.\label{epsW22}
\end{align}
Eventually, $\eps^\delta_{\mathrm{W},4}(t,\varphi)$ satisfies the same estimate as $\eps^\delta_{\mathrm{W},2}(t,\varphi)$ with $\partial_\xi\varphi$ instead of $\varphi$ in the right-hand side of \eqref{epsW22}.
\label{lem:consisSTO}\end{lemma}

\textbf{Proof of Lemma~\ref{lem:consisSTO}.} Define
\begin{align*}
\eps^\delta_{\mathrm{W},1}(t,\varphi)&=\int_0^t\int_{\T^N}\int_{\R\times\R} \left[g_{k,\delta}(x,\xi)-g_k(x,\xi)\right]\varphi(x,\zeta)d\mu^\delta_{x,s,t}(\xi,\zeta) dx d\beta_k(s),
\end{align*}
and let $\eps^\delta_{\mathrm{W},2}(t,\varphi)$ be equal to
$$
\int_0^t\int_{\T^N}\left[\int_{\R\times\R} g_k(x,\xi)\varphi(x,\zeta)d\mu^\delta_{x,s,t}(\xi,\zeta)-\int_\R g_k(x,\xi)\varphi(x,\xi)d\nu^\delta_{x,s}(\xi) \right]dx d\beta_k(s).
$$
Then \eqref{consistSTONoise} is satisfied. Note that
$
n\mapsto\eps^\delta_{\mathrm{W},1}(t_n,\varphi)
$
is a $(\mathcal{F}_{t_n})$-martingale. By Doob's Inequality, Jensen's Inequality (note that $\mu^\delta_{x,s,t}(\R\times\R)\leq 1$) and \eqref{D1num0}, we deduce
\begin{align*}
&\E\left[\sup_{0\leq n<N_T}|\eps^\delta_{\mathrm{W},1}(t_n,\varphi)|^2\right]\\
\leq&4\E\int_0^{t_{N_T}}\left|\int_{\T^N}\int_{\R\times\R} \left[g_{k,\delta}(x,\xi)-g_k(x,\xi)\right]\varphi(x,\zeta)d\mu^\delta_{x,s,t_{N_T}}(\xi,\zeta)\right|^2 dx ds\\
\leq &4\E\int_0^{t_{N_T}}\int_{\T^N}\int_{\R\times\R} \left|g_{k,\delta}(x,\xi)-g_k(x,\xi)\right|^2d\mu^\delta_{x,s,t}(\xi,\zeta) dx ds\|\varphi\|_{L^\infty_{x,\xi}}^2\\
\leq & 4 D_1 T \|\varphi\|_{L^\infty_{x,\xi}}^2 h^2.
\end{align*}
Besides, we see, using It\={o}'s Isometry, and \eqref{D0}, \eqref{D0num}, that
\begin{align*}
&\E\left[\sup_{t\in[t_n,t_{n+1})}|\eps^\delta_{\mathrm{W},1}(t,\varphi)-\eps^\delta_{\mathrm{W},1}(t_n,\varphi)|^2\right]\\
\leq &\E\int_{t_n}^{t_{n+1}}\left|\int_{\T^N}\int_{\R\times\R} \left[g_{k,\delta}(x,v_\delta(s,x))-g_k(x,v_\delta(s,x))\right]\varphi(x,v^\sharp_\delta(s,x))\right|^2 dx ds\\
\leq &2D_0 \|\varphi\|_{L^\infty_{x,\xi}}^2\left[\sup_{0\leq n<N_T}\Delta t_n\right].
\end{align*}
Similarly, we have
$$
\E\left[\sup_{t\in[0,T]}|\eps^\delta_{\mathrm{W},2}(t,\varphi)|^2\right]
\leq 2D_0 \|\varphi\|_{L^\infty_{x,\xi}}^2\left[\sup_{0\leq n<{N_T}}\Delta t_n\right]+
\E\left[\sup_{0\leq n<{N_T}}|\eps^\delta_{\mathrm{W},2}(t_n,\varphi)|^2\right].
$$
Using Doob's Inequality, we obtain
$$
\E\left[\sup_{t\in[0,T]}|\eps^\delta_{\mathrm{W},2}(t,\varphi)|^2\right]
\leq 2D_0 \|\varphi\|_{L^\infty_{x,\xi}}^2\left[\sup_{0\leq n<{N_T}}\Delta t_n\right]+4\E|\eps^\delta_{\mathrm{W},2}(t_{N_T},\varphi)|^2.
$$
By It\={o}'s Formula, $\E|\eps^\delta_{\mathrm{W},2}(t_{N_T},\varphi)|^2$ is bounded from above by
\begin{equation}\label{epsSTO20}
\E\int_0^{t_{N_T}}\int_{\T^N}\sum_{k\geq 1}\left|\int_{\R\times\R} g_k(x,\xi)\varphi(x,\zeta)d\mu^\delta_{x,s,t_{N_T}}(\xi,\zeta)-\int_\R g_k(x,\xi)\varphi(x,\xi)d\nu^\delta_{x,s}(\xi) \right|^2 dx ds. 
\end{equation}
We have, for $t\in [0,T)$, $t\in[t_n,t_{n+1})$, $n<{N_T}$, and $\psi\in C_b(\R\times\R)$,
\begin{align*}
\<\mu^\delta_{x,t,t_{N_T}},\psi\>
=&\<\nu_{x,t}\otimes\nu_{x,t},\psi\>
-\frac{t-t_n}{\Delta t_n}\left[\psi(v^\sharp_\delta(x,t),v^\sharp_\delta(x,t))-\psi(v_\delta(x,t),v^\sharp_\delta(x,t))\right]\\
&-\frac{t_{n+1}-t}{\Delta t_n}\left[\psi(v_\delta(x,t),v_\delta(x,t))-\psi(v_\delta(x,t),v^\sharp_\delta(x,t))\right].
\end{align*}
We estimate therefore \eqref{epsSTO20} by the two terms
\begin{equation}\label{epsSTO23}
2\E\sum_{0\leq n<{N_T}}\int_{t_n}^{t_{n+1}}\int_{\T^N}\sum_{k\geq 1}\left|g_k(x,v^\sharp_\delta(x,t))-g_k(x,v_\delta(x,t))\right|^2|\varphi(x,\xi)|^2 dx dt,
\end{equation}
and
\begin{equation}\label{epsSTO24}
2\E\sum_{0\leq n<{N_T}}\int_{t_n}^{t_{n+1}}\int_{\T^N}\sum_{k\geq 1}\left|\varphi(x,v^\sharp_\delta(x,t))-\varphi(x,v_\delta(x,t))\right|^2 |g_k(x,v_\delta(x,t))|^2 dx dt.
\end{equation}
Note that \eqref{D1} gives, for all $\eta>0$, and $\overline{v}$, $v\in\R$,
\begin{equation}\label{epsSTO25}
\sum_{k\geq 1}\left|g_k(x,\overline{v})-g_k(x,v)\right|^2\leq D_1|\overline{v}-v|\leq D_1\left(\eta+\frac{1}{\eta}|\overline{v}-v|^2\right).
\end{equation}
In virtue of \eqref{epsSTO25}, we can bound \eqref{epsSTO23} by
\begin{equation*}
2D_1\|\varphi\|_{L^\infty_{x,v}}^2\left[\eta+\frac{1}{\eta}\E\int_0^T\|v^\sharp_\delta(t)-v_\delta(t)\|_{L^2({\T^N})^2}^2 dt\right].
\end{equation*}
Using \eqref{timeBVvbarv} and taking $\eta=\left[\sup_{0\leq n<{N_T}}\Delta t_n\right]^{1/2}$, we deduce that \eqref{epsSTO23}  is bounded by
 \begin{equation}\label{epsSTO26}
2D_1\|\varphi\|_{L^\infty_{x,\xi}}^2\left[\frac{1}{\theta}\|v_\delta(0)\|_{L^2({\T^N})}^2+\frac{3 D_0 T}{\theta}\right]\left[\sup_{0\leq n<{N_T}}\Delta t_n\right]^{1/2}.
\end{equation}
An estimate on \eqref{epsSTO24} is obtained as follows: \eqref{epsSTO24} is bounded by
$$
2D_0\|\partial_\xi\varphi\|_{L^\infty_{x,\xi}}^2\E\int_0^T\|v^\sharp_\delta(t)-v_\delta(t)\|_{L^2({\T^N})^2}^2 dt.
$$
Using \eqref{timeBVvbarv} gives an estimate on \eqref{epsSTO24} from above by
 \begin{equation}\label{epsSTO27}
2D_0\|\partial_\xi\varphi\|_{L^\infty_{x,\xi}}^2\left[\frac{1}{\theta}\|v_\delta(0)\|_{L^2({\T^N})}^2+\frac{3 D_0 T}{\theta}\right]\left[\sup_{0\leq n<{N_T}}\Delta t_n\right].
\end{equation}
Next, we denote by $\eps^\delta_{\mathrm{W},3}(t,\varphi)$ and $\eps^\delta_{\mathrm{W},4}(t,\varphi)$ the error terms
\begin{align*}
&\int_0^t\int_{\T^N}\int_{\R\times\R} \left[\GG^2_{k,\delta}(x,\xi)-\GG^2_k(x,\xi)\right]\partial_\xi\varphi(x,\zeta)d\mu^\delta_{x,s,t}(\xi,\zeta) dx ds,\\
&\int_0^t\int_{\T^N}\left[\int_{\R\times\R} \GG^2_k(x,\xi)\partial_\xi\varphi(x,\zeta)d\mu^\delta_{x,s,t}(\xi,\zeta)-\int_\R \GG^2_k(x,\xi)\partial_\xi\varphi(x,\xi)d\nu^\delta_{x,s}(\xi) \right]dx ds.
\end{align*}
We have, for $x\in K$, $\eta>0$,
\begin{align*}
|\GG^2_{k,\delta}(x,\xi)-\GG^2_k(x,\xi)|&=\left|\sum_{k\geq 1}(g_{k,K}(\xi)-g_k(x,\xi))(g_{k,K}(\xi)+g_k(x,\xi))\right|\\
&\leq\frac{1}{2\eta}\sum_{k\geq 1}|g_{k,K}(\xi)-g_k(x,\xi)|^2+\eta\sum_{k\geq 1}|g_{k,K}(\xi)|^2+|g_k(x,\xi)|^2. 
\end{align*}
Using \eqref{D1num0}, \eqref{D0} and \eqref{D0num} and taking $\eta=h$, we see that
\begin{equation}\label{D0D1}
|\GG^2_{k,\delta}(x,\xi)-\GG^2_k(x,\xi)|\leq (D_0+D_1)h.
\end{equation}
This is sufficient to obtain \eqref{epsW33} and the last statement of the lemma (estimate on $\eps^\delta_{\mathrm{W},4}(t,\varphi)$). \qed

\subsection{Conclusion}\label{sec:appcl}

To conclude, let us set
$$
\eps_\delta(\varphi)=\sum_{j=1}^4 \eps^\delta_{\mathrm{W},j}(t,\varphi)-\eps^\delta_{\mathrm{space},0}(t,\varphi)-\eps^\delta_{\mathrm{space},1}(t,\varphi).
$$ 
Then the approximate kinetic equation \eqref{ExactCorrection} follows from the discrete kinetic equation \eqref{AppCorrection} and from the consistency estimates \eqref{consistSPACE}-\eqref{consistSTONoise}-\eqref{consistSTOCorrection}. Since $\|v_\delta(0)\|_{L^2({\T^N})}\leq \|u_0\|_{L^2({\T^N})}$ (the projection \eqref{FVIC} onto piecewise-constant functions is an orthogonal projection in $L^2({\T^N})$), it follows from the error estimates \eqref{eps00SPACE}, \eqref{eps11SPACE}, \eqref{epsW11}, \eqref{epsW33}, \eqref{epsW22} and from the CFL condition \eqref{CFLstrongU} that 
$$
\E\left[\sup_{t\in[0,T]}|\eps^\delta(t,\varphi)|^2\right]\leq C(\varphi)|\delta|^{1/2},
$$
where $C(\varphi)$ is a constant that depends on $\|u_0\|_{L^2({\T^N})}$, on $D_0$, $D_1$, ${L_A}$, on the parameter $\theta$ in \eqref{CFLstrongU}, on $T$, on $|\Lambda|$, where $\Lambda$ is the support of $\varphi$, and on the norms $\|\partial_{x_i}^{j_i}\partial_\xi^k\varphi\|_{L^\infty({\T^N}\times\R)}$ with $j_i+k\leq 2$.

\section{Convergence}\label{sec:Convergence}

To apply Theorem~\ref{th:pathcv} on the basis of Proposition~\ref{prop:consis}, we need to establish some additional estimates on the numerical Young measure $\nu^\delta$ and on the numerical random measure $m^\delta$. This is done in Section~\ref{sec:addestim}. We conclude to the convergence of the Finite Volume method in Section~\ref{sec:Convergence}, Theorem~\ref{th:mainthm}.

\subsection{Additional estimates}\label{sec:addestim}

\subsubsection{Tightness of $(\nu^\delta)$}\label{sec:tightnuh}

\begin{lemma}[Tightness of $(\nu^\delta)$] Let $u_0\in L^\infty({\T^N})$, $T>0$ and $\delta\in \mathfrak{d}_T$. Assume that \eqref{multiplicativeNoise}, \eqref{D0plus}, \eqref{AALip} and \eqref{CFLstrongU} are satisfied. Let $(v_\delta(t))$ be the numerical unknown defined by \eqref{FVscheme}-\eqref{FVIC}-\eqref{defvh} and let $\nu^\delta$ be defined by \eqref{defnuh}. Let $p\in[1,+\infty)$. We have
\begin{equation}\label{eq:tightnuh}
\E\left(\sup_{t\in[0,T]}\int_{{\T^N}}\int_{\R}(1+|\xi|^p) d\nu^\delta_{x,t}(\xi) dx\right)\leq C_p,
\end{equation}
where $C_p$ is a constant depending on $D_0$, $p$, $T$ and $\|u_0\|_{L^\infty({\T^N})}$ only.
\label{lem:tightnuh}\end{lemma}

\textbf{Proof of Lemma~\ref{lem:tightnuh}.} It is sufficient to do the proof for $p\in 2\N^*$ since
$
1+|\xi|^p\leq 2(1+|\xi|^q)
$
for all $\xi\in\R$ if $q\geq p$.
Note that 
$$
\int_{{\T^N}}\int_{\R}|\xi|^p d\nu^\delta_{x,t}(\xi)dx=
\frac{t-t_n}{\Delta t_n}\|v^\sharp_\delta(t)\|_{L^p({\T^N})}^p+\frac{t_{n+1}-t}{\Delta t_n}\|v_\delta(t)\|_{L^p({\T^N})}^p,
$$
for $t\in[t_n,t_{n+1})$. Recall also that $v^\sharp_\delta$ is defined by \eqref{def:barvni}. Let 
$$
\varphi_p(\xi)=p\xi^{p-1}=\partial_\xi\xi^p\quad\xi\in\R.
$$ 
We multiply Equation~\eqref{KiFVhalf} by $\varphi_p(\xi)$ and sum the result over $K$, $\xi$. We obtain then, thanks to \eqref{conservativesyma},
\begin{equation}\label{Lphalf}
\|v_\delta^\flat(t_{n+1})\|_{L^p({\T^N})}^p+p(p-1)\Delta t_n\sum_{K\in\mathcal{T}/\Z^N}|K|\int_\R \xi^{p-2}m^n_K(\xi)d\xi= \|v_\delta(t_{n})\|_{L^p({\T^N})}^p.
\end{equation}
In particular, we have the $L^p$ estimate
\begin{equation}\label{Lphalfleq}
\|v_\delta^\flat(t_{n+1})\|_{L^p({\T^N})}^p\leq \|v_\delta(t_{n})\|_{L^p({\T^N})}^p.
\end{equation}
Let us now estimate the increase of $L^p$-norm due to the stochastic evolution. By It\={o}'s Formula and \eqref{def:barvni}, we have 
\begin{multline*}
|v^{n+1}_K|^p\\
=|v^{n+1/2}_K|^p+p\int_{t_n}^{t_{n+1}}v_K^\sharp(t)^{p-1} g_{k,K}(v^n_K)d\beta_k(t)+\frac12 p(p-1)\int_{t_n}^{t_{n+1}}v_K^\sharp(t)^{p-2} \GG^2_K(v^n_K)dt,
\end{multline*}
and thus
\begin{multline}\label{Lphalf2}
\|v_\delta(t_{n+1})\|_{L^p({\T^N})}^p= \|v_\delta^\flat(t_{n+1})\|_{L^p({\T^N})}^p
+p\int_{t_n}^{t_{n+1}}\<v^\sharp_\delta(t)^{p-1},\gamma_k^n\>_{L^2({\T^N})}d\beta_k(t)\\
+\frac12 p(p-1)\int_{t_n}^{t_{n+1}}\<v^\sharp_\delta(t)^{p-2},\Gamma^n\>_{L^2({\T^N})}dt,
\end{multline}
where 
$$
\gamma_k^n(x)=g_{k,K}(v^n_K),\quad \Gamma^n(x)=\GG^2_K(v^n_K),\quad x\in K.
$$
Using \eqref{Lphalfleq} and induction, we obtain
\begin{equation}\label{LpMartingale}
\|v_\delta(T)\|_{L^p({\T^N})}^p\leq \|v_\delta(0)\|_{L^p({\T^N})}^p+M_{N_T}+B_{N_T},
\end{equation}
where $(M_N)$ is the martingale
\begin{equation*}
M_N=p\sum_{n=0}^{N-1}\int_{t_n}^{t_{n+1}}\<v^\sharp_\delta(t)^{p-1},\gamma_k^n\>_{L^2({\T^N})}d\beta_k(t)
\end{equation*}
and
\begin{equation*}
B_N=\frac12 p(p-1)\sum_{n=0}^{N-1}\int_{t_n}^{t_{n+1}}\<v^\sharp_\delta(t)^{p-2},\Gamma^n\>_{L^2({\T^N})}dt.
\end{equation*}
Note that the argument $\<v^\sharp_\delta(t)^{p-2},\Gamma^n\>_{L^2({\T^N})}$ in $B_N$ is non-negative since $\Gamma^n\geq 0$ and $p-2\in 2\N$. Consequently, $\E\sup_{0\leq n\leq {N_T}} B_n=\E B_{N_T}$ and we deduce the following bound
\begin{equation}
\E\sup_{0\leq n\leq {N_T}} B_n
\leq \frac12 p(p-1)D_0\sum_{n=0}^{N_T-1}\int_{t_n}^{t_{n+1}}\E\|v^\sharp_\delta(t)\|_{L^{p-2}({\T^N})}^{p-2} dt.\label{BND0}
\end{equation}
We have used \eqref{D0num} to obtain \eqref{BND0}.
If $p=2$, then $\E\|v^\sharp_\delta(t)\|_{L^{p-2}({\T^N})}^{p-2}=1$ and is therefore bounded. To estimate $\E\|v^\sharp_\delta(t)\|_{L^{p-2}({\T^N})}^{p-2}$ when $p\geq 4$, note that $v^\sharp_\delta(t)=v_\delta^\flat(t_{n+1})+z^n_\delta(t)$ for $t\in(t_n,t_{n+1})$, where $z^n_\delta(x,t):=\gamma^n_k(x) (\beta_k(t)-\beta_k(t_n))$ is, conditionally to $\mathcal{F}_n$, a Gaussian random variable with variance, for $x\in K$,
$$
\E\left[|z^n_K(t)|^2|\mathcal{F}_n\right]=(t-t_n)\GG_K^2(v^n_K)\leq D_0\Delta t_n
$$
by \eqref{D0num}. In particular, we have the bound
\begin{align*}
\E\|z^n_\delta(t)\|_{L^{p-2}({\T^N})}^{p-2}
&=\sum_{K\in\mathcal{T}/\Z^N} |K|\E\left(\E\left[|z^n_K(t)|^{p-2}|\mathcal{F}_n\right]\right)\\
&=C(p)\sum_{K\in\mathcal{T}/\Z^N} |K|\E\left(\E\left[|z^n_K(t)|^2|\mathcal{F}_n\right]\right)^{(p-2)/2}
\leq C(p)(D_0\Delta t_n)^{(p-2)/2},
\end{align*}
where $C(p)$ is a constant depending on $p$. It follows, using \eqref{Lphalfleq}, that we have the estimate
\begin{align*}
\E\|v^\sharp_\delta(t)\|_{L^{p-2}({\T^N})}^{p-2}&\leq C(p,D_0)\left(1+\E\|v_\delta^\flat(t_{n+1})\|_{L^{p-2}({\T^N})}^{p-2}\right)\\
&\leq C(p,D_0)\left(1+\E\|v_\delta(t_n)\|_{L^{p-2}({\T^N})}^{p-2}\right),
\end{align*}
where $C(p,D_0)$ is a constant depending on $p$ and $D_0$. In particular, we have
\begin{equation}\label{barvh}
\sup_{0\leq n<{N_T}}\sup_{t\in(t_n,t_{n+1})}\E\|v^\sharp_\delta(t)\|_{L^{p-2}({\T^N})}^{p-2}\leq C(p,D_0)\left(1+\sup_{0\leq n<{N_T}}\E\|v_\delta(t_n)\|_{L^{p-2}({\T^N})}^{p-2}\right).
\end{equation}
By \eqref{BND0}, we conclude that
\begin{equation}
\E\sup_{0\leq n\leq {N_T}} B_n
\leq C(p,D_0)T\left(1+\sup_{0\leq n<{N_T}}\E\|v_\delta(t_n)\|_{L^{p-2}({\T^N})}^{p-2}\right),\label{BND0p}
\end{equation}
for possibly a different constant $C(p,D_0)$. Let us now turn to the estimate of the quantity $\E\sup_{0\leq n\leq {N_T}} |M_n|$. The martingale $(M_N)$ can be rewritten as a stochastic integral (with an integrand which is a simple function). Consequently, the quadratic variation of $M_{N_T}$ is
\begin{align*}
\<M_{N_T}\>&=p^2\sum_{n=0}^{N_T-1}\int_{t_n}^{t_{n+1}}\sum_k |\<|v^\sharp_\delta(t)|^{p-1},\gamma_k^n\>_{L^2({\T^N})}|^2 dt\nonumber\\
&\leq p^2\sum_{n=0}^{N_T-1}\int_{t_n}^{t_{n+1}}\sum_k \| |v^\sharp_\delta(t)|^{p-1}\|_{L^2({\T^N})}^2\|\gamma_k^n\|_{L^2({\T^N})}^2 dt\nonumber\\
&=p^2\sum_{n=0}^{N_T-1}\int_{t_n}^{t_{n+1}} \|v^\sharp_\delta(t)\|_{L^{2(p-1)}({\T^N})}^{2(p-1)}\|\Gamma^n\|_{L^1({\T^N})} dt\nonumber\\
&\leq p^2D_0 \sum_{n=0}^{N_T-1}\int_{t_n}^{t_{n+1}} \|v^\sharp_\delta(t)\|_{L^{2(p-1)}({\T^N})}^{2(p-1)} dt,
\end{align*}
by \eqref{D0num}. Using \eqref{barvh} (with $2p$ instead of $p$) gives thus
\begin{equation}\label{quadM}
\E\<M_{N_T}\>\leq p^2D_0 T C(2p,D_0)\left(1+\sup_{0\leq n<{N_T}}\E\|v_\delta(t_n)\|_{L^{2p-2}({\T^N})}^{2p-2}\right).
\end{equation} 
By Burkholder - Davis - Gundy's Inequality, there exists a constant $C_{\mathrm{BDG}}$ such that
$$
\E\sup_{0\leq n\leq {N_T}} |M_n| \leq C_{\mathrm{BDG}}\E \<M_{N_T}\>^{1/2}.
$$
By Jensen's Inequality and the estimate \eqref{quadM}, we obtain
\begin{align}
\E\sup_{0\leq n\leq {N_T}} |M_n| &\leq C_{\mathrm{BDG}}(\E \<M_{N_T}\>)^{1/2}\nonumber\\
&\leq C_{\mathrm{BDG}} p (D_0 T C(2p,D_0))^{1/2}\left(1+\sup_{0\leq n<{N_T}}\E\|v_\delta(t_n)\|_{L^{2p-2}({\T^N})}^{2p-2}\right)^{1/2}.\label{MNsup}
\end{align}
We can conclude now. Since $\E M_{N_T}=0$, taking expectation in \eqref{LpMartingale} (where we replace ${N_T}$ by $n$) gives
$$
\E\|v_\delta(t_n)\|_{L^p({\T^N})}^p\leq \|v_\delta(0)\|_{L^p({\T^N})}^p+\E B_n.
$$
Note (see Section~\ref{sec:appcl}) that
$$
\|v_\delta(0)\|_{L^p({\T^N})}\leq\|u_0\|_{L^p({\T^N})}\leq \|u_0\|_{L^\infty({\T^N})}.
$$
By \eqref{BND0p}, this gives 
$$
\sup_{0\leq n<{N_T}}\E\|v_\delta(t_n)\|_{L^p({\T^N})}^p\leq \|u_0\|_{L^\infty({\T^N})}^p+C(p,D_0)T\left(1+\sup_{0\leq n<{N_T}}\E\|v_\delta(t_n)\|_{L^{p-2}({\T^N})}^{p-2}\right).
$$
By iteration on $p\in 2\N^*$, we deduce, for every such $p$, that
\begin{equation}\label{supE}
\sup_{0\leq n<{N_T}}\E\|v_\delta(t_n)\|_{L^p({\T^N})}^p\leq C_p,
\end{equation}
where the constant $C_p$ depends on $p$, $D_0$, $T$ and $\|u_0\|_{L^\infty({\T^N})}$. Denote generally by $C_p$ any such constant, possibly different from line to line, depending only on $p$, $D_0$, $T$ and $\|u_0\|_{L^\infty({\T^N})}$. By \eqref{supE} with $2p-2$ instead of $p$, we have $\E\sup_{0\leq n\leq {N_T}} |M_n|\leq C_p$. Then we use \eqref{BND0p} with $p-2$ instead of $p$ to obtain $\E\sup_{0\leq n\leq {N_T}} B_n\leq C_p$. By \eqref{LpMartingale}, we deduce 
\begin{equation}\label{Esup}
\E \sup_{0\leq n< {N_T}} \|v_\delta(t_n)\|_{L^p({\T^N})}^p\leq C_p,
\end{equation}
which concludes the proof of the lemma. \qed

\subsubsection{Tightness of $(m_\delta)$}\label{sec:tightmh}

\begin{lemma}[Tightness of $(m_\delta)$] Let $u_0\in L^\infty({\T^N})$, $T>0$ and $\delta\in \mathfrak{d}_T$. Assume that \eqref{multiplicativeNoise}, \eqref{D0plus}, \eqref{AALip} and \eqref{CFLstrongU} are satisfied. Let $(v_\delta(t))$ be the numerical unknown defined by \eqref{FVscheme}-\eqref{FVIC}-\eqref{defvh} and let $m_\delta$ be defined by \eqref{defmh}. Then, for all $p\geq 1$, we have
\begin{equation}\label{estimmhadd}
\E\left|\iiint_{{\T^N}\times[0,T)\times\R}(1+|\xi|^p) dm_\delta(x,t,\xi)\right|^2\leq C_p,
\end{equation}
where $C_p$ is a constant depending on $D_0$, $p$, $T$ and $\|u_0\|_{L^\infty({\T^N})}$ only.
\label{lem:tightmh}\end{lemma}

\textbf{Proof of Lemma~\ref{lem:tightmh}.} Let $p\in 2\N^*$. By \eqref{Lphalf}, we have 
\begin{equation*}
\sum_{K\in\mathcal{T}/\Z^N}\int_\R \xi^{p-2}m^n_K(\xi)d\xi\leq \|v_\delta(t_{n})\|_{L^p({\T^N})}^p.
\end{equation*}
We multiply this inequality by $\Delta t_n$. Summing over $n\in\{0,\ldots,{N_T}-1\}$, we obtain 
$$
\iiint_{{\T^N}\times[0,T)\times\R}\xi^{p-2} dm_\delta(x,t,\xi)\leq T \sup_{0\leq n< {N_T}} \|v_\delta(t_n)\|_{L^p({\T^N})}^p.
$$
Since $\|v_\delta(t_n)\|_{L^p({\T^N})}^{2p}\leq \|v_\delta(t_n)\|_{L^{2p}({\T^N})}^{2p}$, we have
$$
\E\left|\iiint_{{\T^N}\times[0,T)\times\R}\xi^{p-2} dm_\delta(x,t,\xi)\right|^2\leq
T^2 \E\sup_{0\leq n< {N_T}} \|v_\delta(t_n)\|_{L^{2p}({\T^N})}^{2p}\leq T^2 C_p
$$
by \eqref{Esup}. Taking first $p=2$, then $p\in 2\N^*$ arbitrary, we obtain \eqref{estimmhadd}. \qed

\subsection{Convergence}\label{subsec:convergence}

We may now apply the theorem~\ref{th:pathcv}, to obtain the following results.

\begin{theorem} Let $u_0\in L^\infty({\T^N})$, $T>0$. Assume that the hypotheses \eqref{D1}, \eqref{multiplicativeNoise}, \eqref{D0plus}, 
are satisfied. Then there exists a unique solution $u$ to \eqref{stoSCL} with initial datum $u_0$, in the sense of Definition~\ref{defkineticsol}. Besides, for all $1\leq p<+\infty$, almost-surely, $u\in C([0,T];L^p(\T^N))$. 
\label{th:mainthmExist}\end{theorem}

\begin{theorem} Let $u_0\in L^\infty({\T^N})$, $T>0$. Assume that the hypotheses \eqref{D1}, \eqref{multiplicativeNoise}, \eqref{D0plus}, \eqref{AALip}, \eqref{alphaK}, \eqref{alphapK} and \eqref{CFLstrongU} are satisfied. Let $u$ be the solution to \eqref{stoSCL} with initial datum $u_0$ and let $v_\delta$ be the solution to the Finite Volume scheme \eqref{FVscheme}-\eqref{FVIC}-\eqref{Halfflux}-\eqref{defX}. Then we have the convergence
\begin{equation}\label{CVfinal}
\lim_{\delta\to0}\E\|v_\delta-u\|_{L^p(\T^N\times(0,T))}^p=0,
\end{equation}
for all $p\in[1,\infty)$.
\label{th:mainthm}\end{theorem}

\begin{remark}\label{rk:rkcvlaw} If the $X^{n+1}_k$ are merely i.i.d. random variables with normalized centred normal law $\mathcal{N}(0,1)$, then $(v_\delta)$ is converging to $u$ in $L^p(\T^N\times(0,T))$ \textit{in law} when $\delta\to0$. Indeed, the identity \eqref{defX} is only satisfied in law now, hence $v_\delta$ has the same law as the function $\tilde v_\delta$ defined by \eqref{FVscheme}-\eqref{FVIC}-\eqref{Halfflux}, with $X^{n+1}_k$ replaced by the right-hand side of \eqref{defX}. We apply the conclusion of Theorem~\ref{th:mainthm} to $\tilde v_\delta$. As a corollary, we obtain the convergence in law of $(\tilde v_\delta)$ to $u$ in $L^p(\T^N\times(0,T))$. A slightly different manner of expressing the same thing is to notice that, when the discrete increments $(X^{n+1}_k)$ are some given normal law $\mathcal{N}(0,1)$, then we can construct a set of Brownian motions $\tilde\beta_k(t)$ such that
\begin{equation}\label{defXXX}
X^{n+1}_k=\frac{\tilde\beta_k(t_{n+1})-\tilde\beta_k(t_n)}{(\Delta t_n)^{1/2}}.
\end{equation}
Indeed, without loss of generality, we can restrict ourselves to the case $\Delta t_n=1$ in \eqref{defXXX} and use the
L{\'e}vy-Ciesielski construction of the Brownian motion, \cite[Section~3.2]{SchillingPartzsch2014} on $[0,1]$ as follows: we define (\textit{cf.} \cite[Formula~(3.1)]{SchillingPartzsch2014} )
$$
G_0=X^{n+1}_k,\; G_1=X^{n+1}_1,\ldots, G_{k-1}=X^{n+1}_{k-1},\; G_k=X^{n+1}_{k+1},\; G_{k+1}=X^{n+1}_{k+2},\ldots
$$
and we set
$$
\tilde{\beta}(t)=\sum_{p=0}^\infty G_p \<\mathbf{1}_{[0,t]},H_p\>,
$$
where the $H_p$'s are the Haar basis of $L^2(0,1)$. Then \eqref{defXXX} follows from the fact that
$$
\int_0^1 H_p(t) dt=\<H_p,H_0\>=\delta_{p0}.
$$
\end{remark}

\begin{remark}\label{rk:2th} Theorem~\ref{th:mainthmExist} has already been proved in \cite{DebusscheVovelle10} (see also Section~5 in \cite{DottiVovelle16a}) under less restrictive hypotheses (having a compactly supported noise is unnecessary). We give the statement together with Theorem~\ref{th:mainthm}, however, to emphasize the fact that the convergence of the Finite Volume method, in the framework which we use, give both the existence-uniqueness of the solution to the limit continuous problem, and the convergence of the numerical method to this solution. It is not necessary to provide the existence of the solution $u$ to \eqref{stoSCL} by an external means.
\end{remark} 

\textbf{Proof of Theorem~\ref{th:mainthmExist} and Theorem~\ref{th:mainthm}.} Let us first prove the theorem~\ref{th:mainthm}. We take the existence of $u$, solution to \eqref{stoSCL} with initial datum $u_0$ for granted. By Proposition~\ref{prop:consis}, Lemma~\ref{lem:tightnuh} and Lemma~\ref{lem:tightmh}, we may apply Theorem~\ref{th:pathcv} to $f_\delta$: we obtain \eqref{CVfinal} with $z_\delta$ instead of $v_\delta$, where
$$
z_\delta(x,t):=\int_\R\xi d\nu^\delta_{x,t}(\xi)=\frac{t-t_n}{\Delta t_n}\overline{v}_\delta(x,t)+\frac{t_{n+1}-t}{\Delta t_n}v_\delta(x,t),
$$
for $t\in[t_n,t_{n+1}]$. By \eqref{timeBVvbarv}, we have the estimate 
$$
\E\int_0^T\|z_\delta-v_\delta\|^2_{L^2(\T^N)}dt=\mathcal{O}(|\delta|)
$$
on the difference between $z_\delta$ and $v_\delta$. This gives \eqref{CVfinal} for $p\leq 2$. 
If $p>2$, we use the inequality
\begin{align}
\E\|v_\delta-u\|_{L^p(\T^N\times(0,T))}^p&\leq \|v_\delta-u\|_{L^2(\Omega\times\T^N\times(0,T))}\|v_\delta-u\|_{L^{2(p-1)}(\Omega\times\T^N\times(0,T))}^{p-1}\nonumber\\
&\leq \frac{1}{\eta}\E\|v_\delta-u\|_{L^2(\T^N\times(0,T))}^2+\eta \E\|v_\delta-u\|_{L^{2(p-1)}(\T^N\times(0,T))}^{2(p-1)},\label{prefinal}
\end{align}
where $\eta$ is a positive parameter. Thanks to the uniform bounds \eqref{eq:integrabilityu}-\eqref{Esup}, we can choose $\eta$ independent on $\delta$ to have the second term in \eqref{prefinal} smaller than an arbitrary threshold. 
By the convergence result for $p=2$ the first term in \eqref{prefinal} is then also small for $\delta$ close to $0$. This concludes the proof of theorem~\ref{th:mainthm}. To prove \ref{th:mainthmExist}, we just need to construct an approximation scheme satisfying \eqref{AALip}, \eqref{alphaK}, \eqref{alphapK} and \eqref{CFLstrongU} and to compute $v_\delta$ by \eqref{FVscheme}-\eqref{FVIC}-\eqref{Halfflux}-\eqref{defX}. We can use a cartesian grid for this purpose: let $h_m=\frac{1}{m}$, where $m\in\N^*$. Let ${\mathcal{T}/\Z^N}$ be the set of open hypercubes of length $h_m$ obtained by translates of the original hypercube $(0,h_m)^N$ by vectors $h_m x$, $x$ having components in $\{0,l\dots,m-1\}$. Then \eqref{alphaK}, \eqref{alphapK} are satisfied with $\alpha_N=2^{-N}$ since an hypercube has $2^N$ sides. We can choose the Godunov numerical fluxes, defined as follows:
$$
A_{K\to L}(v,w)=\begin{cases}
\displaystyle |K|L|\min_{v\leq\xi\leq w}A(\xi)\cdot n_{K,L} & \mbox{if }v\leq w,\\
\displaystyle |K|L|\max_{w\leq\xi\leq v}A(\xi)\cdot n_{K,L} & \mbox{if }w\leq v.
\end{cases}
$$
These fluxes $A_{K\to L}$ are monotone (non-increasing in the first variable, non-increasing in the second variable). They satisfy the hypotheses of regularity, consistency \eqref{AALip} and \eqref{consistency} respectively with $L_A=\mathrm{Lip}(A)$. The conservative symmetry property \eqref{conservativesym} is also satisfied. At last, to ensure \eqref{CFLstrongU}, we just need to take a uniform time step $\Delta t$ like
$$
\Delta t=\frac12 \frac{\alpha_N^2}{2 {L_A}}\, h_m.
$$
We have then \eqref{CFLstrongU} with $\theta=\frac12$. At this stage, Proposition~\ref{prop:consis} provides a sequence of approximate solutions $(f_m)$. Thanks to the uniform bounds established in Section~\ref{sec:addestim}, we can apply Theorem~\ref{th:pathcv}: this gives the existence of a unique solution $u$ to \eqref{stoSCL} with initial datum $u_0$. By Corollary~3.3 in \cite{DottiVovelle16a}, we have $u\in C([0,T];L^p(\T^N))$ almost-surely. \qed


\def\ocirc#1{\ifmmode\setbox0=\hbox{$#1$}\dimen0=\ht0 \advance\dimen0
  by1pt\rlap{\hbox to\wd0{\hss\raise\dimen0
  \hbox{\hskip.2em$\scriptscriptstyle\circ$}\hss}}#1\else {\accent"17 #1}\fi}
  \def\cprime{$'$}

\end{document}